\documentclass{article}
\usepackage{bm,amsmath,amsthm,url} 
\usepackage[dvips]{graphicx}
\RequirePackage[colorlinks,citecolor=blue,urlcolor=blue]{hyperref}
\theoremstyle{plain}
\newtheorem{theorem}{Theorem}
\newtheorem{corollary}[theorem]{Corollary} 
\newtheorem*{conjecture}{Conjecture} 
\theoremstyle{remark} 
\newtheorem*{remark}{Remark} 
\def\P{{\rm P}} 
\def\E{{\rm E}} 
\def\Var{{\rm Var}} 
\def\Cov{{\rm Cov}} 
\def\T{{\sf T}} 
\def\eps{\varepsilon} 
\allowdisplaybreaks
\makeatletter 
\renewcommand\@biblabel[1]{}
\makeatother
\begin{document}

\title{Limit theorems for Parrondo's paradox}

\author{S. N. Ethier\\
\begin{small}University of Utah\end{small}\\
\begin{small}Department of Mathematics\end{small}\\
\begin{small}155 S. 1400 E., JWB 233\end{small}\\
\begin{small}Salt Lake City, UT 84112, USA\end{small}\\
\begin{small}e-mail: \url{ethier@math.utah.edu}\end{small}
\and
Jiyeon Lee\thanks{Supported by the Yeungnam University research grants in 2008.}\\
\begin{small}Yeungnam University\end{small}\\
\begin{small}Department of Statistics\end{small}\\
\begin{small}214-1 Daedong, Kyeongsan\end{small}\\
\begin{small}Kyeongbuk 712-749, South Korea\end{small}\\
\begin{small}e-mail: \url{leejy@yu.ac.kr}\end{small}}
\date{}
\maketitle

\begin{abstract}
That there exist two losing games that can be combined, either by random mixture or by nonrandom alternation, to form a winning game is known as Parrondo's paradox.  We establish a strong law of large numbers and a central limit theorem for the Parrondo player's sequence of profits, both in a one-parameter family of capital-dependent games and in a two-parameter family of history-dependent games, with the potentially winning game being either a random mixture or a nonrandom pattern of the two losing games.  We derive formulas for the mean and variance parameters of the central limit theorem in nearly all such scenarios; formulas for the mean permit an analysis of when the Parrondo effect is present.
\bigskip

 \medskip

  \medskip

\noindent Submitted to EJP on February 18, 2009, final version accepted August 3, 2009.
\end{abstract}

\newpage

\section{Introduction}

The original Parrondo (1996) games can be described as follows: Let $p:={1\over2}-\eps$ and 
\begin{equation}\label{profit}
p_0:={1\over10}-\eps,\qquad p_1:={3\over4}-\eps,
\end{equation}
where $\eps>0$ is a small bias parameter (less than 1/10, of course).  In game $A$, the player tosses a $p$-coin (i.e., $p$ is the probability of heads).  In game $B$, if the player's current capital is divisible by 3, he tosses a $p_0$-coin, otherwise he tosses a $p_1$-coin.  (Assume initial capital is 0 for simplicity.) In both games, the player wins one unit with heads and loses one unit with tails.

It can be shown that games $A$ and $B$ are both losing games, regardless of $\eps$, whereas the random mixture $C:={1\over2} A+{1\over2} B$ (toss a fair coin to determine which game to play) is a winning game for $\eps$ sufficiently small.  Furthermore, certain nonrandom patterns, including $AAB$, $ABB$, and $AABB$ but excluding $AB$, are winning as well, again for $\eps$ sufficiently small.  These are examples of \textit{Parrondo's paradox}.

The terms ``losing'' and ``winning'' are meant in an asymptotic sense.  More precisely, assume that the game (or mixture or pattern of games) is repeated ad infinitum.  Let $S_n$ be the player's cumulative profit after $n$ games for $n\ge1$.  A game (or mixture or pattern of games) is \textit{losing} if $\lim_{n\to\infty}S_n=-\infty$ a.s., it is \textit{winning} if $\lim_{n\to\infty}S_n=\infty$ a.s.,  and it is \textit{fair} if $-\infty=\liminf_{n\to\infty}S_n<\limsup_{n\to\infty}S_n=\infty$ a.s.  These definitions are due in this context to Key, K\l osek, and Abbott (2006).

Because the games were introduced to help explain the so-called flashing Brownian ratchet (Ajdari and Prost 1992), much of the work on this topic has appeared in the physics literature.  Survey articles include Harmer and Abbott (2002), Parrondo and Din\'is (2004), Epstein (2007), and Abbott (2009).

Game $B$ can be described as \textit{capital-dependent} because the coin choice depends on current capital.  An alternative game $B$, called \textit{history-dependent}, was introduced by Parrondo, Harmer, and Abbott (2000):  Let
\begin{equation}\label{history}
p_0:={9\over10}-\eps,\quad p_1=p_2:={1\over4}-\eps,\quad p_3:={7\over10}-\eps,
\end{equation}
where $\eps>0$ is a small bias parameter.  Game $A$ is as before.  In game $B$, the player tosses a $p_0$-coin (resp., a $p_1$-coin, a $p_2$-coin, a $p_3$-coin) if his two previous results are loss-loss (resp., loss-win, win-loss, win-win).  He wins one unit with heads and loses one unit with tails.

The conclusions for the history-dependent games are the same as for the capital-dependent ones, except that the pattern $AB$ need not be excluded.

Pyke (2003) proved a strong law of large numbers for the Parrondo player's sequence of profits in the capital-dependent setting.  In the present paper we generalize his result and obtain a central limit theorem as well.  We formulate a stochastic model general enough to include both the capital-dependent and the history-dependent games.  We also treat separately the case in which the potentially winning game is a random mixture of the two losing games (game $A$ is played with probability $\gamma$, and game $B$ is played with probability $1-\gamma$) and the case in which the potentially winning game (or, more precisely, pattern of games) is a nonrandom pattern of the two losing games, specifically the pattern $[r,s]$, denoting $r$ plays of game $A$ followed by $s$ plays of game $B$.  Finally, we replace (\ref{profit}) by 
$$
p_0:={\rho^2\over1+\rho^2}-\eps,\qquad p_1:={1\over1+\rho}-\eps,
$$
where $\rho>0$; (\ref{profit}) is the special case $\rho=1/3$.
We also replace (\ref{history}) by 
$$
p_0:={1\over 1+\kappa}-\eps,\quad p_1=p_2:={\lambda\over 1+\lambda}-\eps,\quad p_3:=1-{\lambda\over 1+\kappa}-\eps,
$$
where $\kappa>0$, $\lambda>0$, and $\lambda<1+\kappa$; (\ref{history}) is the special case $\kappa=1/9$ and $\lambda=1/3$.  The reasons for these parametrizations are explained in Sections 3 and 4.

Section 2 formulates our model and derives an SLLN and a CLT.  Section 3 specializes to the capital-dependent games and their random mixtures, showing that the Parrondo effect is present whenever $\rho\in(0,1)$ and $\gamma\in(0,1)$.  Section 4 specializes to the history-dependent games and their random mixtures, showing that the Parrondo effect is present whenever either $\kappa<\lambda<1$ or $\kappa>\lambda>1$ and $\gamma\in(0,1)$.  Section 5 treats the nonrandom patterns $[r,s]$ and derives an SLLN and a CLT.  Section 6 specializes to the capital-dependent games, showing that the Parrondo effect is present whenever $\rho\in(0,1)$ and $r,s\ge1$ with one exception: $r=s=1$.  Section 7 specializes to the history-dependent games.  Here we expect that the Parrondo effect is present whenever either $\kappa<\lambda<1$ or $\kappa>\lambda>1$ and $r,s\ge1$ (without exception), but although we can prove it for certain specific values of $\kappa$ and $\lambda$ (such as $\kappa=1/9$ and $\lambda=1/3$), we cannot prove it in general.  Finally, Section 8 addresses the question of why Parrondo's paradox holds.

In nearly all cases we obtain formulas for the mean and variance parameters of the CLT.  These parameters can be interpreted as the asymptotic mean per game played and the asymptotic variance per game played of the player's cumulative profit.  Of course, the pattern $[r,s]$ comprises $r+s$ games.  

Some of the algebra required in what follows is rather formidable, so we have used \textit{Mathematica 6} where necessary.  Our \texttt{.nb} files are available upon request.

\section{A general formulation of Parrondo's games}

In some formulations of Parrondo's games, the player's cumulative profit $S_n$ after $n$ games is described by some type of random walk $\{S_n\}_{n\ge1}$, and then a Markov chain $\{X_n\}_{n\ge0}$ is defined in terms of $\{S_n\}_{n\ge1}$; for example, $X_n\equiv \xi_0+S_n\;({\rm mod}\;3)$ in the capital-dependent games, where $\xi_0$ denotes initial capital.  However, it is more logical to introduce the Markov chain $\{X_n\}_{n\ge0}$ first and then define the random walk $\{S_n\}_{n\ge1}$ in terms of $\{X_n\}_{n\ge0}$.

Consider an irreducible aperiodic Markov chain $\{X_n\}_{n\ge0}$ with finite state space $\Sigma$.  It evolves according to the one-step transition matrix\footnote{In the physics literature the one-step transition matrix is often written in transposed form, that is, with column sums equal to 1.  We do not follow that convention here.  More precisely, here $P_{ij}:=\P(X_n=j\mid X_{n-1}=i)$.} ${\bm P}=(P_{ij})_{i,j\in\Sigma}$.  Let us denote its unique stationary distribution by ${\bm \pi}=(\pi_i)_{i\in \Sigma}$.  Let $w:\Sigma\times\Sigma\mapsto {\bf R}$ be an arbitrary function, which we will sometimes write as a matrix ${\bm W}:=(w(i,j))_{i,j\in\Sigma}$ and refer to as the \textit{payoff matrix}. Finally, define the sequences $\{\xi_n\}_{n\ge1}$ and $\{S_n\}_{n\ge1}$ by
\begin{equation}\label{xi_n}
\xi_n:=w(X_{n-1},X_n),\qquad n\ge1,
\end{equation}
and
\begin{equation}\label{S_n}
S_n:=\xi_1+\cdots+\xi_n,\qquad n\ge1.
\end{equation}

For example, let $\Sigma:=\{0,1,2\}$, put $X_0:=\xi_0\;({\rm mod\ }3)$, $\xi_0$ being initial capital, and let the payoff matrix be given by\footnote{Coincidentally, this is the payoff matrix for the classic game \textit{stone-scissors-paper}.  However, Parrondo's games, as originally formulated, are games of chance, not games of strategy, and so are outside the purview of game theory (in the sense of von Neumann).}
\begin{equation}\label{W-profit}
\setlength{\arraycolsep}{1.5mm}
{\bm W}:=\left(\begin{array}{rrr}
0&1&-1\\
-1&0&1\\
1&-1&0
\end{array}\right).
\end{equation}
With the role of ${\bm P}$ played by
$$
\setlength{\arraycolsep}{1.5mm}
{\bm P}_B:=\left(\begin{array}{ccc}
0&p_0&1-p_0\\
1-p_1&0&p_1\\
p_1&1-p_1&0
\end{array}\right),
$$
where $p_0$ and $p_1$ are as in (\ref{profit}),
$S_n$ represents the player's profit after $n$ games when playing the capital-dependent game $B$ repeatedly.  With the role of ${\bm P}$ played by
$$
\setlength{\arraycolsep}{1.5mm}
{\bm P}_A:=\left(\begin{array}{ccc}
0&p&1-p\\
1-p&0&p\\
p&1-p&0
\end{array}\right),
$$
where $p:={1\over2}-\eps$,
$S_n$ represents the player's profit after $n$ games when playing game $A$ repeatedly.  Finally, with the role of ${\bm P}$ played by ${\bm P}_C:=\gamma{\bm P}_A+(1-\gamma){\bm P}_B$, where $0<\gamma<1$, $S_n$ represents the player's profit after $n$ games when playing the mixed game $C:=\gamma A+(1-\gamma)B$ repeatedly.  In summary, all three capital-dependent games are described by the same stochastic model with a suitable choice of parameters.  

Similarly, the history-dependent games fit into the same framework, as do the ``primary'' Parrondo games of Cleuren and Van den Broeck (2004).

Thus, our initial goal is to analyze the asymptotic behavior of $S_n$ under the conditions of the second paragraph of this section.  We begin by assuming that $X_0$ has distribution ${\bm \pi}$, so that $\{X_n\}_{n\ge0}$ and hence $\{\xi_n\}_{n\ge1}$ are stationary sequences, although we will weaken this assumption later.  

We claim that the conditions of the stationary, strong mixing central limit theorem apply to $\{\xi_n\}_{n\ge1}$.  To say that $\{\xi_n\}_{n\ge1}$ has the \textit{strong mixing property} (or is \textit{strongly mixing}) means that $\alpha(n)\to0$ as $n\to\infty$, where
$$
\alpha(n):=\sup_m\sup_{E\in\sigma(\xi_1,\ldots,\xi_m),\,F\in\sigma(\xi_{m+n},\xi_{m+n+1},\ldots)}|\P(E\cap F)-\P(E)\P(F)|.
$$
A version of the theorem for bounded sequences (Bradley 2007, Theorem 10.3) suffices here.  That version requires that $\sum_{n\ge1} \alpha(n)<\infty$.  In our setting, $\{X_n\}_{n\ge0}$ is strongly mixing with a geometric rate (Billingsley 1995, Example 27.6), hence so is $\{\xi_n\}_{n\ge1}$.  Since $\xi_1,\xi_2,\ldots$ are bounded by $\max|w|$, it follows that
$$
\sigma^2:=\Var(\xi_1)+2\sum_{m=1}^\infty\Cov(\xi_1,\xi_{m+1})
$$
converges absolutely.  For the theorem to apply, it suffices to assume that $\sigma^2>0$.

Let us evaluate the mean and variance parameters of the central limit theorem.  First,
\begin{equation}\label{mu}
\mu:=\E[\xi_1]=\sum_i \P(X_0=i)\E[w(X_0,X_1)\mid X_0=i]
=\sum_{i,j} \pi_i P_{ij}w(i,j)
\end{equation}
and
$$
\Var(\xi_1)=\E[\xi_1^2]-(\E[\xi_1])^2
=\sum_{i,j} \pi_i P_{ij}w(i,j)^2-\bigg(\sum_{i,j} \pi_i P_{ij}w(i,j)\bigg)^2.
$$
To evaluate $\Cov(\xi_1,\xi_{m+1})$, we first find
\begin{eqnarray*}
\E[\xi_1\xi_{m+1}]&=&\sum_i\pi_i\E[w(X_0,X_1)\E[w(X_m,X_{m+1})\mid X_0,X_1]\mid X_0=i]\\
&=&\sum_{i,j}\pi_i P_{ij}w(i,j) \E[w(X_m,X_{m+1})\mid X_1=j]\\
&=&\sum_{i,j,k,l}\pi_i P_{ij}w(i,j) ({\bm P}^{m-1})_{jk}P_{kl}w(k,l),
\end{eqnarray*}
from which it follows that
$$
\Cov(\xi_1,\xi_{m+1})=\sum_{i,j,k,l}\pi_i P_{ij}w(i,j)[({\bm P}^{m-1})_{jk}-\pi_k]P_{kl}w(k,l).
$$
We conclude that
$$
\sum_{m=1}^\infty \Cov(\xi_1,\xi_{m+1})
=\sum_{i,j,k,l}\pi_i P_{ij}w(i,j)(z_{jk}-\pi_k)P_{kl}w(k,l),
$$
where ${\bm Z}=(z_{ij})$ is the \textit{fundamental matrix} associated with ${\bm P}$ (Kemeny and Snell 1960, p.\ 75).  

In more detail, we let ${\bm \Pi}$ denote the square matrix each of whose rows is ${\bm \pi}$, and we find that
$$
\sum_{m=1}^\infty ({\bm P}^{m-1}-{\bm \Pi})={\bm I}-{\bm \Pi}+\sum_{n=1}^\infty ({\bm P}^{n}-{\bm \Pi})={\bm Z}-{\bm \Pi},
$$
where
\begin{equation}\label{Z}
{\bm Z}:={\bm I}+\sum_{n=1}^\infty ({\bm P}^n-{\bm \Pi})=({\bm I}-({\bm P}-{\bm \Pi}))^{-1};
\end{equation}
for the last equality, see Kemeny and Snell (loc.\ cit.).  Therefore,
\begin{eqnarray}\label{sigma2}
\sigma^2&=&\sum_{i,j} \pi_i P_{ij}w(i,j)^2-\bigg(\sum_{i,j} \pi_i P_{ij}w(i,j)\bigg)^2\nonumber \\
&&\quad{}+2\sum_{i,j,k,l}\pi_i P_{ij}w(i,j)(z_{jk}-\pi_k)P_{kl}w(k,l).
\end{eqnarray}
A referee has pointed out that these formulas can be written more concisely using matrix notation.  Denote by $\bm P'$ (resp., $\bm P''$) the matrix whose $(i,j)$th entry is $P_{ij}w(i,j)$ (resp., $P_{ij}w(i,j)^2$), and let $\bm 1:=(1,1,\ldots,1)^\T$.  Then
\begin{equation}\label{mu,sigma2}
\mu=\bm\pi\bm P'\bm 1\quad{\rm and}\quad\sigma^2=\bm\pi\bm P''\bm 1
-(\bm\pi\bm P'\bm 1)^2+2\bm\pi\bm P'(\bm Z-\bm\Pi)\bm P'\bm 1.
\end{equation}

If $\sigma^2>0$, then (Bradley 2007, Proposition 8.3)
$$
\lim_{n\to\infty}n^{-1}\Var(S_n)=\sigma^2
$$
and the central limit theorem applies, that is, $(S_n-n\mu)/\sqrt{n\sigma^2}\to_d N(0,1)$.  If we strengthen the assumption that $\sum_{n\ge1}\alpha(n)<\infty$ by assuming that $\alpha(n)=O(n^{-(1+\delta)})$, where $0<\delta<1$, we have (Bradley 2007, proof of Lemma 10.4)
$$
\E[(S_n-n\mu)^4]=O(n^{3-\delta}),
$$
hence by the Borel--Cantelli lemma it follows that $S_n/n\to\mu$ a.s.  In other words, the strong law of large numbers applies.  

Finally, we claim that, if $\mu=0$ and $\sigma^2>0$, then
\begin{equation}\label{liminf-limsup}
-\infty=\liminf_{n\to\infty}S_n<\limsup_{n\to\infty}S_n=\infty\;\;{\rm a.s.}
\end{equation}
Indeed, $\{\xi_n\}_{n\ge1}$ is stationary and strongly mixing, hence its future tail $\sigma$-field is trivial (Bradley 2007, p.\ 60), in the sense that every event has probability 0 or 1.  It follows that $\P(\liminf_{n\to\infty}S_n=-\infty)$ is 0 or 1.  Since $\mu=0$ and $\sigma^2>0$, we can invoke the central limit theorem to conclude that this probability is 1.  Similarly, we get $\P(\limsup_{n\to\infty}S_n=\infty)=1$. 

Each of these derivations required that the sequence $\{\xi_n\}_{n\ge1}$ be stationary, an assumption that holds if $X_0$ has distribution ${\bm \pi}$, but in fact the distribution of $X_0$ can be arbitrary, and $\{\xi_n\}_{n\ge1}$ need not be stationary. 

\begin{theorem}
Let $\mu$ and $\sigma^2$ be as in (\ref{mu}) and (\ref{sigma2}).
Under the assumptions of the second paragraph of this section, but with the distribution of $X_0$ arbitrary,
\begin{equation}\label{SLLN}
\lim_{n\to\infty}n^{-1}\E[S_n]=\mu\quad\text{and}\quad{S_n\over n}\to \mu\;\;{\rm a.s.}
\end{equation}
and, if $\sigma^2>0$, 
\begin{equation}\label{CLT}
\lim_{n\to\infty}n^{-1}\Var(S_n)=\sigma^2\quad\text{and}\quad{S_n-n\mu\over\sqrt{n\sigma^2}}\to_d N(0,1).
\end{equation}
If $\mu=0$ and $\sigma^2>0$, then (\ref{liminf-limsup}) holds.
\end{theorem}

\begin{remark}
Assume that $\sigma^2>0$.  It follows that, if $S_n$ is the player's cumulative profit after $n$ games for each $n\ge1$, then the game (or mixture or pattern of games) is losing if $\mu<0$, winning if $\mu>0$, and fair if $\mu=0$.  (See Section 1 for the definitions of these three terms.)
\end{remark}

\begin{proof} 
It will suffice to treat the case $X_0=i_0\in\Sigma$, and then use this case to prove the general case.
Let $\{X_n\}_{n\ge0}$ be a Markov chain in $\Sigma$ with one-step transition matrix ${\bm P}$ and initial distribution ${\bm \pi}$, so that $\{\xi_n\}_{n\ge1}$ is stationary, as above.  
Let $N:=\min\{n\ge0: X_n=i_0\}$, and define
$$
\hat X_n:=X_{N+n},\qquad n\ge0.
$$
Then $\{\hat X_n\}_{n\ge0}$ is a Markov chain in $\Sigma$ with one-step transition matrix ${\bm P}$ and initial state $\hat X_0=i_0$. We can define $\{\hat\xi_n\}_{n\ge1}$ and $\{\hat S_n\}_{n\ge1}$ in terms of it by analogy with (\ref{xi_n}) and (\ref{S_n}). If $n\ge N$, then
\begin{eqnarray*}
\hat S_n-S_n&=&\hat\xi_1+\cdots+\hat\xi_n-(\xi_1+\cdots+\xi_n)\\
&=&\hat\xi_1+\cdots+\hat\xi_{n-N}+\hat\xi_{n-N+1}+\cdots+\hat\xi_n\\
&&\qquad{}-(\xi_1+\cdots+\xi_N+\xi_{N+1}+\cdots+\xi_n)\\
&=&\hat\xi_{n-N+1}+\cdots+\hat\xi_n-(\xi_1+\cdots+\xi_N),
\end{eqnarray*}
and $\hat S_n-S_n$ is bounded by $2N\max |w|$.  Thus, if we divide by $n$ or $\sqrt{n\sigma^2}$, the result tends to 0 a.s.\ as $n\to\infty$.  We get the SLLN and the CLT with $\hat S_n$ in place of $S_n$, via this coupling of the two sequences.  We also get the last conclusion in a similar way.  The first equation in (\ref{SLLN}) follows from the second using bounded convergence.  Finally, the random variable $N$ has finite moments (Durrett 1996, Chapter 5, Exercise 2.5).  The first equation in (\ref{CLT}) therefore follows from our coupling.  
\end{proof}

A well-known (Bradley 2007, pp.\ 36--37) nontrivial example for which $\sigma^2=0$ is the case in which $\Sigma\subset{\bf R}$ and $w(i,j)=j-i$ for all $i,j\in\Sigma$.  Then $S_n=X_n-X_0$ (a telescoping sum) for all $n\ge1$, hence $\mu=0$ and $\sigma^2=0$ by Theorem 1.

However, it may be of interest to confirm these conclusions using only the mean, variance, and covariance formulas above.  We calculate
$$
\mu:=\E[\xi_1]=\sum_{i,j}\pi_i P_{ij}(j-i)=\sum_j\pi_j j-\sum_i\pi_i i=0,
$$
\begin{equation}\label{eqn1}
\Var(\xi_1)=\sum_{i,j}\pi_i P_{ij}(j-i)^2-\mu^2
=2\sum_{i}\pi_i i^2-2\sum_{i,j}\pi_i P_{ij}ij,
\end{equation}
and
\begin{eqnarray}\label{eqn2}
\sum_{m=1}^\infty\Cov(\xi_1,\xi_{m+1})
&=&\sum_{i,j,k,l}\pi_i P_{ij}(j-i)z_{jk}P_{kl}(l-k)\nonumber\\
&=&\sum_{i,j,k,l}\pi_i P_{ij}z_{jk}P_{kl}jl-\sum_{i,j,k,l}\pi_i P_{ij}z_{jk}P_{kl}jk\nonumber\\
&&\quad{}-\sum_{i,j,k,l}\pi_i P_{ij}z_{jk}P_{kl}il+\sum_{i,j,k,l}\pi_i P_{ij}z_{jk}P_{kl}ik\nonumber\\
&=&\sum_{j,k,l}\pi_j z_{jk}P_{kl}jl-\sum_{j,k}\pi_j z_{jk}jk\nonumber\\
&&\quad{}-\sum_{i,j,k,l}\pi_i P_{ij}z_{jk}P_{kl}il+\sum_{i,j,k}\pi_i P_{ij}z_{jk}ik.
\end{eqnarray}
Since $\bm Z=(\bm I-(\bm P-\bm\Pi))^{-1}$, we can multiply $(\bm I-\bm P+\bm \Pi)\bm Z=\bm I$ on the right by $\bm P$ and use $\bm\Pi\bm Z\bm P=\bm\Pi\bm P=\bm\Pi$ to get $\bm P\bm Z\bm P=\bm Z\bm P+\bm\Pi-\bm P$. 
This implies that
$$
\sum_{i,j,k,l}\pi_i P_{ij}z_{jk}P_{kl}il=\sum_{i,k,l}\pi_i z_{ik}P_{kl}il+\sum_{i,l}\pi_i\pi_l il-\sum_{i,l}\pi_i P_{il}il.
$$
From the fact that $\bm P\bm Z=\bm Z+\bm\Pi-\bm I$, we also obtain
$$
\sum_{i,j,k}\pi_i P_{ij}z_{jk}ik=\sum_{i,k}\pi_i z_{ik}ik+\sum_{i,k}\pi_i\pi_k ik-\sum_{i}\pi_i i^2.
$$
Substituting in (\ref{eqn2}) gives
$$
\sum_{m=1}^\infty\Cov(\xi_1,\xi_{m+1})
=\sum_{i,l}\pi_i P_{il}il-\sum_{i}\pi_i i^2,
$$
and this, together with (\ref{eqn1}), implies that $\sigma^2=0$.

\section{Mixtures of capital-dependent games}

The Markov chain underlying the capital-dependent Parrondo games has state space $\Sigma=\{0,1,2\}$ and one-step transition matrix of the form
\begin{equation}\label{P-profit}
\setlength{\arraycolsep}{1.5mm}
{\bm P}:=\left(\begin{array}{ccc}
0&p_0&1-p_0\\
1-p_1&0&p_1\\
p_2&1-p_2&0
\end{array}\right),
\end{equation}
where $p_0,p_1,p_2\in(0,1)$.  It is irreducible and aperiodic.  The payoff matrix ${\bm W}$ is as in (\ref{W-profit}).

It will be convenient below to define $q_i:=1-p_i$ for $i=0,1,2$.
Now, the unique stationary distribution ${\bm\pi}=(\pi_0,\pi_1,\pi_2)$ has the form
$$
\pi_0=(1-p_1q_2)/d,\quad
\pi_1=(1-p_2q_0)/d,\quad
\pi_2=(1-p_0q_1)/d,
$$
where $d:=2+p_0p_1p_2+q_0q_1q_2$.
Further, the fundamental matrix ${\bm Z}=(z_{ij})$ is easy to evaluate (e.g., $z_{00}=\pi_0+[\pi_1(1+p_1)+\pi_2(1+q_2)]/d$).

We conclude that $\{\xi_n\}_{n\ge1}$ satisfies the SLLN with 
$$
\mu=\sum_{i=0}^2\pi_i(p_i-q_i)
$$
and the CLT with the same $\mu$ and with
$$
\sigma^2=1-\mu^2+2\sum_{i=0}^2\sum_{k=0}^2\pi_i [p_i (z_{[i+1],k}-\pi_k)
-q_i (z_{[i-1],k}-\pi_k)](p_k-q_k),
$$
where $[j]\in\{0,1,2\}$ satisfies $j\equiv[j]$ (mod 3), at least if $\sigma^2>0$.

We now apply these results to the capital-dependent Parrondo games. 
Although actually much simpler, game $A$ fits into this framework with one-step transition matrix ${\bm P}_A$ defined by (\ref{P-profit}) with
$$
p_0=p_1=p_2:={1\over2}-\eps,
$$
where $\eps>0$ is a small bias parameter.
In game $B$, it is typically assumed that, ignoring the bias parameter, the one-step transition matrix $\bm P_B$ is defined by (\ref{P-profit}) with
$$
p_1=p_2\quad{\rm and}\quad \mu=0.
$$
These two constraints determine a one-parameter family of probabilities given by 
\begin{equation}\label{profit-param}
p_1=p_2={1\over 1+\sqrt{p_0/(1-p_0)}}.
\end{equation}
To eliminate the square root, we reparametrize the probabilities in terms of $\rho>0$.  Restoring the bias parameter, game $B$ has one-step transition matrix ${\bm P}_B$ defined by (\ref{P-profit}) with
\begin{equation}\label{probs-profit}
p_0:={\rho^2\over1+\rho^2}-\eps,\qquad p_1=p_2:={1\over 1+\rho}-\eps,
\end{equation}
which includes (\ref{profit}) when $\rho=1/3$.
Finally, game $C:=\gamma A+(1-\gamma)B$ is a mixture $(0<\gamma<1)$ of the two games, hence has one-step transition matrix ${\bm P}_C:=\gamma{\bm P}_A+(1-\gamma){\bm P}_B$, which can also be defined by (\ref{P-profit}) with
$$
p_0:=\gamma\,{1\over2}+(1-\gamma){\rho^2\over1+\rho^2}-\eps,\qquad p_1=p_2:=\gamma\,{1\over2}+(1-\gamma){1\over 1+\rho}-\eps.
$$

Let us denote the mean $\mu$ for game $A$ by $\mu_A(\eps)$, to emphasize the game as well as its dependence on $\eps$.  Similarly, we denote the variance $\sigma^2$ for game $A$ by $\sigma_A^2(\eps)$.  Analogous notation applies to games $B$ and $C$.  We obtain, for game $A$, $\mu_A(\eps)=-2\eps$ and $\sigma_A^2(\eps)=1-4\eps^2$; for game $B$,
$$
\mu_B(\eps)=-{3(1+2\rho+6\rho^2+2\rho^3+\rho^4)\over2(1+\rho+\rho^2)^2}\,\eps +O(\eps^2)
$$
and 
$$
\sigma_B^2(\eps)=\left({3\rho\over1+\rho+\rho^2}\right)^2+O(\eps);
$$
and for game $C$,
$$
\mu_C(\eps)={3\gamma(1-\gamma)(2-\gamma)(1-\rho)^3(1+\rho)\over
8(1+\rho+\rho^2)^2+\gamma(2-\gamma)(1-\rho)^2(1+4\rho+\rho^2)}+O(\eps)
$$
and 
\begin{eqnarray*}
\sigma_C^2(\eps)&=&1-\mu_C(0)^2-32(1-\gamma)^2(1-\rho^3)^2[16(1+\rho+\rho^2)^2(1+4\rho+\rho^2)\\
&&\quad{}+8\gamma(1-\rho)^2(1+4\rho+\rho^2)^2+24\gamma^2(1-\rho)^2(1-\rho-6\rho^2-\rho^3+\rho^4)\\
&&\quad{}-\gamma^3(4-\gamma)(1-\rho)^4(7+16\rho+7\rho^2)]\\
&&\qquad\qquad\;\;{}/[8(1+\rho+\rho^2)^2+\gamma(2-\gamma)(1-\rho)^2(1+4\rho+\rho^2)]^3+O(\eps).
\end{eqnarray*}
One can check that $\sigma_C^2(0)<1$ for all $\rho\ne1$ and $\gamma\in(0,1)$.

The formula for $\sigma_B^2(0)$ was found by Percus and Percus (2002) in a different form.  We prefer the form given here because it tells us immediately that game $B$ has smaller variance than game $A$ for each $\rho\ne1$, provided $\eps$ is sufficiently small.  With $\rho=1/3$ and $\eps=1/200$, Harmer and Abbott (2002, Fig.~5) inferred this from a simulation.

The formula for $\mu_C(0)$ was obtained by Berresford and Rockett (2003) in a different form.  We prefer the form given here because it makes the following conclusion transparent.

\begin{theorem}[Pyke 2003]
Let $\rho>0$ and let games $A$ and $B$ be as above but with the bias parameter absent. If $\gamma\in(0,1)$ and $C:=\gamma A+(1-\gamma)B$, then $\mu_C(0)>0$ for all $\rho\in(0,1)$, $\mu_C(0)=0$ for $\rho=1$, and $\mu_C(0)<0$ for all $\rho>1$.
\end{theorem}

Assuming (\ref{probs-profit}) with $\eps=0$, the condition $\rho<1$ is equivalent to $p_0<{1\over2}$.

Clearly, the Parrondo effect appears (with the bias parameter present) if and only if $\mu_C(0)>0$.  A reverse Parrondo effect, in which two winning games combine to lose, appears (with a \textit{negative} bias parameter present) if and only if $\mu_C(0)<0$.

\begin{corollary}[Pyke 2003]
Let games $A$ and $B$ be as above (with the bias parameter present).
If $\rho\in(0,1)$ and $\gamma\in(0,1)$, then there exists $\eps_0>0$, depending on $\rho$ and $\gamma$, such that Parrondo's paradox holds for games $A$, $B$, and $C:=\gamma A+(1-\gamma)B$, that is, $\mu_A(\eps)<0$, $\mu_B(\eps)<0$, and $\mu_C(\eps)>0$, whenever $0<\eps<\eps_0$.
\end{corollary}

The theorem and corollary are special cases of results of Pyke.  In his form\-ulation, the modulo 3 condition in the definition of game $B$ is replaced by a modulo $m$ condition, where $m\ge3$, and game $A$ is replaced by a game analogous to game $B$ but with a different parameter $\rho_0$.  Pyke's condition is equivalent to $0<\rho<\rho_0\le1$.  We have assumed $m=3$ and $\rho_0=1$.

We would like to point out a useful property of game $B$.  We assume $\eps=0$ and we temporarily denote $\{X_n\}_{n\ge0}$, $\{\xi_n\}_{n\ge1}$, and $\{S_n\}_{n\ge1}$ by $\{X_n(\rho)\}_{n\ge0}$, $\{\xi_n(\rho)\}_{n\ge1}$, and $\{S_n(\rho)\}_{n\ge1}$ to emphasize their dependence on $\rho$.  Similarly, we temporarily denote $\mu_B(0)$ and $\sigma_B^2(0)$ by $\mu_B(\rho,0)$ and $\sigma_B^2(\rho,0)$.  Replacing $\rho$ by $1/\rho$ has the effect of changing the win probabilities $p_0=\rho^2/(1+\rho^2)$ and $p_1=1/(1+\rho)$ to the loss probabilities $1-p_0$ and $1-p_1$, and vice versa.  Therefore, given $\rho\in(0,1)$, we expect that
$$
\xi_n(1/\rho)=-\xi_n(\rho), \quad  S_n(1/\rho)=-S_n(\rho),\qquad n\ge1,
$$
a property that is in fact realized by coupling the Markov chains $\{X_n(\rho)\}_{n\ge0}$ and $\{X_n(1/\rho)\}_{n\ge0}$ so that $X_n(1/\rho)\equiv-X_n(\rho)\;({\rm mod}\;3)$ for all $n\ge1$.  It follows that
\begin{equation}\label{symmetry}
\mu_B(1/\rho,0)=-\mu_B(\rho,0)\quad{\rm and}\quad\sigma_B^2(1/\rho,0)=\sigma_B^2(\rho,0).
\end{equation}
The same argument gives (\ref{symmetry}) for game $C$.  The reader may have noticed that the formulas derived above for $\mu_B(0)$ and $\sigma_B^2(0)$, as well as those for $\mu_C(0)$ and $\sigma_C^2(0)$, satisfy (\ref{symmetry}).

When $\rho=1/3$, the mean and variance formulas simplify to
$$
\mu_B(\eps)=-{294\over169}\,\eps+O(\eps^2)\quad{\rm and}\quad\sigma_B^2(\eps)={81\over169}+O(\eps)
$$
and, if $\gamma={1\over2}$ as well, 
$$
\mu_C(\eps)={18\over709}+O(\eps)\quad{\rm and}\quad\sigma_C^2(\eps)=\frac{311313105}{356400829}+O(\eps).
$$

\section{Mixtures of history-dependent games}

The Markov chain underlying the history-dependent Parrondo games has state space $\Sigma=\{0,1,2,3\}$ and one-step transition matrix of the form
\begin{equation}\label{P-history}
\setlength{\arraycolsep}{1.5mm}
{\bm P}:=\left(\begin{array}{cccc}
1-p_0&p_0&0&0\\
0&0&1-p_1&p_1\\
1-p_2&p_2&0&0\\
0&0&1-p_3&p_3
\end{array}\right),
\end{equation}
where $p_0,p_1,p_2,p_3\in(0,1)$.  Think of the states of $\Sigma$ in binary form: 00, 01, 10, 11.  They represent, respectively, loss-loss, loss-win, win-loss, and win-win for the results of the two preceding games, with the second-listed result being the more recent one.  The Markov chain is irreducible and aperiodic.  The payoff matrix ${\bm W}$ is given by
\begin{equation*}\label{W-history}
\setlength{\arraycolsep}{1.5mm}
{\bm W}:=\left(\begin{array}{rrrr}
-1&1&0&0\\
0&0&-1&1\\
-1&1&0&0\\
0&0&-1&1
\end{array}\right).
\end{equation*}

It will be convenient below to define $q_i:=1-p_i$ for $i=0,1,2,3$.
Now, the unique stationary distribution ${\bm \pi}=(\pi_0,\pi_1,\pi_2,\pi_3)$ has the form
$$
\pi_0=q_2q_3/d,\quad
\pi_1=p_0q_3/d,\quad
\pi_2=p_0q_3/d,\quad
\pi_3=p_0p_1/d,
$$
where $d:=p_0p_1+2p_0q_3+q_2q_3$.
Further, the fundamental matrix $\bm Z=(z_{ij})$ can be evaluated with some effort (e.g., $z_{00}=\pi_0+[\pi_1(p_1+2q_3)+\pi_2(p_1p_2+p_2q_3+q_3)+\pi_3(p_1p_2+p_2q_3+q_2+q_3)]/d$).

We conclude that $\{\xi_n\}_{n\ge1}$ satisfies the SLLN with 
$$
\mu=\sum_{i=0}^3\pi_i(p_i-q_i)
$$
and the CLT with the same $\mu$ and with
$$
\sigma^2=1-\mu^2+2\sum_{i=0}^3\sum_{k=0}^3\pi_i [p_i (z_{[2i+1],k}-\pi_k)
-q_i (z_{[2i],k}-\pi_k)](p_k-q_k),
$$
where $[j]\in\{0,1,2,3\}$ satisfies $j\equiv[j]$ (mod 4), at least if $\sigma^2>0$.

We now apply these results to the history-dependent Parrondo games.  Although actually much simpler, game $A$ fits into this framework with one-step transition matrix ${\bm P}_A$ defined by (\ref{P-history}) with
$$
p_0=p_1=p_2=p_3:={1\over2}-\eps,
$$
where $\eps>0$ is a small bias parameter.  In game $B$, it is typically assumed that, ignoring the bias parameter, the one-step transition matrix $\bm P_B$ is defined by (\ref{P-history}) with
$$
p_1=p_2\quad{\rm and}\quad \mu=0.
$$
These two constraints determine a two-parameter family of probabilities given by 
\begin{equation}\label{history-param}
p_1=p_2\quad{\rm and}\quad p_3=1-{p_0p_1\over1-p_1}.
\end{equation}
We reparametrize the probabilities in terms of $\kappa>0$ and $\lambda>0$ (with $\lambda<1+\kappa$).  Restoring the bias parameter, game $B$ has one-step transition matrix ${\bm P}_B$ defined by (\ref{P-history}) with
\begin{equation}\label{probs-history}
p_0:={1\over 1+\kappa}-\eps,\quad p_1=p_2:={\lambda\over 1+\lambda}-\eps,\quad
p_3:=1-{\lambda\over 1+\kappa}-\eps,
\end{equation}
which includes (\ref{history}) when $\kappa=1/9$ and $\lambda=1/3$.
Finally, game $C:=\gamma A+(1-\gamma)B$ is a mixture ($0<\gamma<1$) of the two games, hence has one-step transition matrix ${\bm P}_C:=\gamma{\bm P}_A+(1-\gamma){\bm P}_B$, which also has the form (\ref{P-history}).

We obtain, for game $A$, $\mu_A(\eps)=-2\eps$ and $\sigma_A^2(\eps)=1-4\eps^2$; for game $B$,
$$
\mu_B(\eps)=-{(1+\kappa)(1+\lambda)\over2\lambda}\,\eps +O(\eps^2)
$$
and 
$$
\sigma_B^2(\eps)={(1+\kappa)(1+\kappa+\kappa\lambda+\kappa\lambda^2)\over \lambda(1+\lambda)(2+\kappa+\lambda)}+O(\eps);
$$
and for game $C$,
\begin{eqnarray*}
\mu_C(\eps)&=&\gamma(1-\gamma)(1+\kappa)(\lambda-\kappa)(1-\lambda)/[2\gamma(2-\gamma)+\gamma(5-\gamma)\kappa\\
&&\quad{}+(8-9\gamma+3\gamma^2)\lambda+\gamma(1+\gamma)\kappa^2+4\kappa\lambda+(1-\gamma)(4-\gamma)\lambda^2\\
&&\quad{}+\gamma(1+\gamma)\kappa^2\lambda+3\gamma(1-\gamma)\kappa\lambda^2]+O(\eps)
\end{eqnarray*}
and 
\begin{eqnarray*}
\sigma_C^2(\eps)&=&1-\mu_C(0)^2+8(1-\gamma)[2-\gamma(1-\kappa)][2\lambda+\gamma(1+\kappa-2\lambda)]\\
&&\quad{}\cdot[2\lambda(2+\kappa+\lambda)^2(1+2\kappa-2\lambda+\kappa^2-3\lambda^2+\kappa^2\lambda-\lambda^3+\kappa^2 \lambda^2)\\
&&\qquad{}+\gamma(2+\kappa+\lambda)(2+5\kappa-11\lambda+4\kappa^2-16 \kappa \lambda+2\lambda^2+\kappa^3\\
&&\qquad{}-3\kappa^2\lambda-2\kappa\lambda^2+18\lambda^3+2 \kappa^3\lambda+10\kappa^2\lambda^2-10\kappa\lambda^3+12\lambda^4\\ 
&&\qquad{}+8\kappa^3\lambda^2-8\kappa^2\lambda^3-13\kappa\lambda^4+3\lambda^5+2\kappa^3\lambda^3-6\kappa^2\lambda^4-2\kappa\lambda^5\\
&&\qquad{}+\kappa^3\lambda^4+\kappa^2 \lambda^5)-\gamma^2(6+14\kappa-14\lambda+9 \kappa^2-24\kappa\lambda-9\lambda^2\\
&&\qquad{}-18\kappa^2\lambda+2\kappa\lambda^2+16 \lambda^3-\kappa^4-12\kappa^3\lambda+33\kappa^2\lambda^2+16\lambda^4\\ 
&&\qquad{}-4\kappa^4 \lambda+16\kappa^3\lambda^2+12\kappa^2\lambda^3-30\kappa\lambda^4+6 \lambda^5-9\kappa^2\lambda^4-16\kappa\lambda^5\\
&&\qquad{}+\lambda^6-4\kappa^4\lambda^3+6\kappa^2 \lambda^5-2\kappa\lambda^6-\kappa^4\lambda^4+4\kappa^3\lambda^5+3 \kappa^2 \lambda^6)\\
&&\qquad{}+2\gamma^3(1-\kappa\lambda)^2(1+\kappa-\lambda-\lambda^2)^2]\\
&&\quad/\{(1+\lambda)[2\gamma(2-\gamma)+\gamma(5-\gamma)\kappa+(8-9\gamma+3\gamma^2)\lambda+\gamma(1+\gamma)\kappa^2\\ 
&&\qquad{}+4\kappa\lambda+(1-\gamma)(4-\gamma)\lambda^2+\gamma(1+\gamma)\kappa^2\lambda+3\gamma(1-\gamma)\kappa\lambda^2]^3\}\\
&&\quad{}+O(\eps).
\end{eqnarray*}
By inspection, we deduce the following conclusion.

\begin{theorem}
Let $\kappa>0$, $\lambda>0$, and $\lambda<1+\kappa$, and
let games $A$ and $B$ be as above but with the bias parameter absent.  If $\gamma\in(0,1)$ and $C:=\gamma A+(1-\gamma)B$, then $\mu_C(0)>0$ whenever $\kappa<\lambda<1$ or $\kappa>\lambda>1$, $\mu_C(0)=0$ whenever $\kappa=\lambda$ or $\lambda=1$, and $\mu_C(0)<0$ whenever $\lambda<\min(\kappa,1)$ or $\lambda>\max(\kappa,1)$.
\end{theorem}

Assuming (\ref{probs-history}) with $\eps=0$, the condition $\kappa<\lambda<1$ is equivalent to
$$
p_0+p_1>1\quad{\rm and}\quad p_1=p_2<{1\over2},
$$
whereas the condition $\kappa>\lambda>1$ is equivalent to
$$
p_0+p_1<1\quad{\rm and}\quad p_1=p_2>{1\over2}.
$$

Again, the Parrondo effect is present if and only if $\mu_C(0)>0$.

\begin{corollary}
Let games $A$ and $B$ be as above (with the bias parameter present).
If $0<\kappa<\lambda<1$ or $\kappa>\lambda>1$, and if $\gamma\in(0,1)$, then there exists $\eps_0>0$, depending on $\kappa$, $\lambda$, and $\gamma$, such that Parrondo's paradox holds for games $A$, $B$, and $C:=\gamma A+(1-\gamma)B$, that is, $\mu_A(\eps)<0$, $\mu_B(\eps)<0$, and $\mu_C(\eps)>0$, whenever $0<\eps<\eps_0$.
\end{corollary}

When $\kappa=1/9$ and $\lambda=1/3$, the mean and variance formulas simplify to
$$
\mu_B(\eps)=-{20\over9}\,\eps+O(\eps^2)\quad{\rm and}\quad\sigma_B^2(\eps)=\frac{235}{198}+O(\eps)
$$
and, if $\gamma={1\over2}$ as well, 
$$
\mu_C(\eps)=\frac{5}{429}+O(\eps)\quad{\rm and}\quad\sigma_C^2(\eps)=\frac{25324040}{26317863}+O(\eps).
$$
Here, in contrast to the capital-dependent games, the variance of game $B$ is greater than that of game $A$.  This conclusion, however, is parameter dependent.

\section{Nonrandom patterns of games}

We also want to consider nonrandom patterns of games of the form $A^rB^s$, in which game $A$ is played $r$ times, then game $B$ is played $s$ times, where $r$ and $s$ are positive integers.  Such a pattern is denoted in the literature by $[r,s]$.  

Associated with the games are one-step transition matrices for Markov chains in a finite state space $\Sigma$, which we will denote by ${\bm P}_A$ and ${\bm P}_B$, and a function $w:\Sigma\times\Sigma\mapsto{\bf R}$.  We assume that ${\bm P}_A$ and ${\bm P}_B$ are irreducible and aperiodic, as are ${\bm P}:={\bm P}_A^r{\bm P}_B^s$, ${\bm P}_A^{r-1}{\bm P}_B^s{\bm P}_A$, \dots, ${\bm P}_B^s{\bm P}_A^r$, \dots, ${\bm P}_B{\bm P}_A^r{\bm P}_B^{s-1}$ (the $r+s$ cyclic permutations of ${\bm P}_A^r{\bm P}_B^s$).  Let us denote the unique stationary distribution associated with ${\bm P}$ by ${\bm \pi}=(\pi_i)_{i\in \Sigma}$.  The driving Markov chain $\{X_n\}_{n\ge0}$ is \textit{time-inhomogeneous}, with one-step transition matrices ${\bm P}_A, {\bm P}_A, \ldots, {\bm P}_A$ ($r$ times), ${\bm P}_B, {\bm P}_B,\ldots, {\bm P}_B$ ($s$ times), ${\bm P}_A, {\bm P}_A, \ldots, {\bm P}_A$ ($r$ times), ${\bm P}_B, {\bm P}_B,\ldots, {\bm P}_B$ ($s$ times), and so on.  Now define $\{\xi_n\}_{n\ge1}$ and $\{S_n\}_{n\ge1}$ by (\ref{xi_n}) and (\ref{S_n}).  What is the asymptotic behavior of $S_n$ as $n\to\infty$?

Let us give $X_0$ distribution ${\bm \pi}$, at least for now.  The time-inhomogeneous Markov chain $\{X_n\}_{n\ge0}$ is of course not stationary, so we define the (time-homogeneous) Markov chain $\{\bm X_n\}_{n\ge0}$ by
\begin{eqnarray}\label{boldXs}
{\bm X}_1&:=&(X_0,X_1,\ldots,X_{r+s}),\nonumber\\
{\bm X}_2&:=&(X_{r+s},X_{r+s+1},\ldots,X_{2(r+s)}),\\
{\bm X}_3&:=&(X_{2(r+s)},X_{2(r+s)+1},\ldots,X_{3(r+s)}),\nonumber\\
&\vdots&\nonumber
\end{eqnarray}
and we note that this is a stationary Markov chain in a subset of $\Sigma^{r+s+1}$. (The overlap between successive vectors is intentional.)  We assume that it is irreducible and aperiodic in that subset, hence it is strongly mixing.  Therefore,
$(\xi_1,\ldots,\xi_{r+s}),(\xi_{r+s+1},\ldots,\xi_{2(r+s)}), (\xi_{2(r+s)+1},\ldots,\xi_{3(r+s)}),\ldots$ is itself a stationary, strongly mixing sequence, and we can apply the stationary, strong mixing CLT to the sequence
\begin{eqnarray}\label{etas}
\eta_1&:=&\xi_1+\cdots+\xi_{r+s},\nonumber\\
\eta_2&:=&\xi_{r+s+1}+\cdots+\xi_{2(r+s)},\\
\eta_3&:=&\xi_{2(r+s)+1}+\cdots+\xi_{3(r+s)},\nonumber\\
&\vdots&\nonumber
\end{eqnarray}
We denote by $\hat\mu$ and $\hat\sigma^2$ the mean and variance parameters for this sequence.

The mean and variance parameters $\mu$ and $\sigma^2$ for the original sequence $\xi_1,\xi_2,\ldots$ can be obtained from these.  Indeed, 
$$
\mu:=\lim_{n\to\infty}{1\over (r+s)n}\E[S_{(r+s)n}]=\lim_{n\to\infty}{1\over (r+s)n}\E[\eta_1+\cdots+\eta_n]={\hat\mu\over r+s},
$$
where $\hat\mu:=\E[\eta_1]$,
and
$$
\sigma^2:=\lim_{n\to\infty}{1\over (r+s)n}\Var(S_{(r+s)n})=\lim_{n\to\infty}{1\over (r+s)n}\Var(\eta_1+\cdots+\eta_n)={\hat\sigma^2\over r+s},
$$
where
$$
\hat\sigma^2:={\Var}(\eta_1)+2\sum_{m=1}^\infty\Cov(\eta_1,\eta_{m+1}).
$$
It remains to evaluate these variances and covariances.  

First,
\begin{equation}\label{mu-periodic}
\mu={1\over r+s}\bigg(\sum_{u=0}^{r-1}\sum_{i,j} ({\bm\pi}{\bm P}_A^u)_i ({\bm P}_A)_{ij}w(i,j)+\sum_{v=0}^{s-1}
\sum_{i,j}({\bm\pi}{\bm P}_A^r{\bm P}_B^v)_i ({\bm P}_B)_{ij}w(i,j)\bigg).
\end{equation}
This formula for $\mu$ is equivalent to one found by Kay and Johnson (2003) in the history-dependent setting.

Next, 
\begin{eqnarray}\label{var}
&&\!\!\!\!\!\!\!\!\Var(\eta_1)\nonumber\\
&=&\sum_{u=0}^{r-1}\Var(\xi_{u+1})+\sum_{v=0}^{s-1}\Var(\xi_{r+v+1})+2\sum_{0\le u<v\le r-1}\Cov(\xi_{u+1},\xi_{v+1})\nonumber\\
&&\;{}+2\sum_{u=0}^{r-1}\sum_{v=0}^{s-1}\Cov(\xi_{u+1},\xi_{r+v+1})
+2\sum_{0\le u<v\le s-1}\Cov(\xi_{r+u+1},\xi_{r+v+1}) \nonumber\\
&=&\sum_{u=0}^{r-1}\bigg[\sum_{i,j}(\bm \pi\bm P_A^u)_i(\bm P_A)_{ij}w(i,j)^2-\bigg(\sum_{i,j}(\bm \pi\bm P_A^u)_i(\bm P_A)_{ij}w(i,j)\bigg)^2\bigg] \nonumber\\
&&\;{}+\sum_{v=0}^{s-1}\bigg[\sum_{i,j}(\bm \pi\bm P_A^r\bm P_B^v)_i(\bm P_B)_{ij}w(i,j)^2-\bigg(\sum_{i,j}(\bm \pi\bm P_A^r\bm P_B^v)_i(\bm P_B)_{ij}w(i,j)\bigg)^2\bigg] \nonumber\\
&&\;{}+2\sum_{0\le u<v\le r-1}\sum_{i,j,k,l}(\bm \pi\bm P_A^u)_i(\bm P_A)_{ij}w(i,j) \nonumber\\
\noalign{\vglue-3mm}
&&\qquad\qquad\qquad\qquad\quad{}\cdot[(\bm P_A^{v-u-1})_{jk}-(\bm \pi\bm P_A^v)_k](\bm P_A)_{kl}w(k,l) \nonumber\\
&&\;{}+2\sum_{u=0}^{r-1}\sum_{v=0}^{s-1}\sum_{i,j,k,l}(\bm \pi\bm P_A^u)_i(\bm P_A)_{ij}w(i,j) \nonumber\\
\noalign{\vglue-3mm}
&&\qquad\qquad\qquad\qquad\quad{}\cdot
[(\bm P_A^{r-u-1}\bm P_B^v)_{jk}-(\bm \pi\bm P_A^r\bm P_B^v)_k](\bm P_B)_{kl}w(k,l) \nonumber\\
&&\;{}+2\sum_{0\le u<v\le s-1}\sum_{i,j,k,l}(\bm \pi\bm P_A^r\bm P_B^u)_i(\bm P_B)_{ij}w(i,j) \nonumber\\
\noalign{\vglue-3mm}
&&\qquad\qquad\qquad\qquad\quad{}\cdot[(\bm P_B^{v-u-1})_{jk}-(\bm \pi\bm P_A^r\bm P_B^v)_k](\bm P_B)_{kl}w(k,l).
\end{eqnarray}

Furthermore,
\begin{eqnarray*}
&&\Cov(\eta_1,\eta_{m+1})\\
&&\;{}=\Cov(\xi_1+\cdots+\xi_{r+s},\xi_{m(r+s)+1}+\cdots+\xi_{(m+1)(r+s)})\\
&&\;{}=\sum_{u=0}^{r-1}\sum_{v=0}^{r-1}\Cov(\xi_{u+1},\xi_{m(r+s)+v+1})+\sum_{u=0}^{r-1}\sum_{v=0}^{s-1}\Cov(\xi_{u+1},\xi_{m(r+s)+r+v+1})\\
&&\quad{}+\sum_{u=0}^{s-1}\sum_{v=0}^{r-1}\Cov(\xi_{r+u+1},\xi_{m(r+s)+v+1})+\sum_{u=0}^{s-1}\sum_{v=0}^{s-1}\Cov(\xi_{r+u+1},\xi_{m(r+s)+r+v+1})\\
&&\;{}=\sum_{u=0}^{r-1}\sum_{v=0}^{r-1}\sum_{i,j,k,l}(\bm \pi\bm P_A^u)_i(\bm P_A)_{ij}w(i,j)\\
\noalign{\vglue-3mm}
&&\qquad\qquad\qquad\quad\;{}\cdot
[(\bm P_A^{r-u-1}\bm P_B^s\bm P^{m-1}\bm P_A^v)_{jk}-(\bm \pi\bm P_A^v)_k](\bm P_A)_{kl}w(k,l)\\
&&\quad{}+\sum_{u=0}^{r-1}\sum_{v=0}^{s-1}\sum_{i,j,k,l}(\bm \pi\bm P_A^u)_i(\bm P_A)_{ij}w(i,j)\\
\noalign{\vglue-3mm}
&&\qquad\qquad\qquad\quad\;{}\cdot[(\bm P_A^{r-u-1}\bm P_B^s\bm P^{m-1}\bm P_A^r\bm P_B^v)_{jk}-(\bm \pi\bm P_A^r\bm P_B^v)_k](\bm P_B)_{kl}w(k,l)\\
&&\quad{}+\sum_{u=0}^{s-1}\sum_{v=0}^{r-1}\sum_{i,j,k,l}(\bm \pi\bm P_A^r\bm P_B^u)_i(\bm P_B)_{ij}w(i,j)\\
\noalign{\vglue-3mm}
&&\qquad\qquad\qquad\quad\;{}\cdot[(\bm P_B^{s-u-1}\bm P^{m-1}\bm P_A^v)_{jk}-(\bm \pi\bm P_A^v)_k](\bm P_A)_{kl}w(k,l)\\
&&\quad{}+\sum_{u=0}^{s-1}\sum_{v=0}^{s-1}\sum_{i,j,k,l}(\bm \pi\bm P_A^r\bm P_B^u)_i(\bm P_B)_{ij}w(i,j)\\
\noalign{\vglue-3mm}
&&\qquad\qquad\qquad\quad\;{}\cdot[(\bm P_B^{s-u-1}\bm P^{m-1}\bm P_A^r\bm P_B^v)_{jk}-(\bm \pi\bm P_A^r\bm P_B^v)_k](\bm P_B)_{kl}w(k,l).
\end{eqnarray*}
Consider the factor $(\bm P_A^{r-u-1}\bm P_B^s\bm P^{m-1}\bm P_A^v)_{jk}-(\bm \pi\bm P_A^v)_k$ in the first sum on the right, for example.  With $\bm \Pi$ denoting the square matrix each of whose rows is $\bm \pi$, we can rewrite this as 
\begin{eqnarray*}
&&(\bm P_A^{r-u-1}\bm P_B^s\bm P^{m-1}\bm P_A^v)_{jk}-(\bm P_A^{r-u-1}\bm P_B^s\bm \Pi\bm P_A^v)_{jk}\\
&&\qquad{}=[\bm P_A^{r-u-1}\bm P_B^s(\bm P^{m-1}-\bm\Pi)\bm P_A^v]_{jk}
\end{eqnarray*}
Thus, summing over $m\ge1$, we get 
\begin{eqnarray}\label{covsum}
&&\sum_{m=1}^\infty\Cov(\eta_1,\eta_{m+1})\nonumber\\
&&\quad{}=\sum_{u=0}^{r-1}\sum_{v=0}^{r-1}\sum_{i,j,k,l}(\bm\pi\bm P_A^u)_i(\bm P_A)_{ij}w(i,j) \nonumber\\
\noalign{\vglue-3mm}
&&\qquad\qquad\qquad\qquad\quad{}\cdot[\bm P_A^{r-u-1}\bm P_B^s(\bm Z-\bm \Pi)\bm P_A^v]_{jk}(\bm P_A)_{kl}w(k,l) \nonumber\\
&&\qquad{}+\sum_{u=0}^{r-1}\sum_{v=0}^{s-1}\sum_{i,j,k,l}(\bm\pi\bm P_A^u)_i(\bm P_A)_{ij}w(i,j) \nonumber\\
\noalign{\vglue-3mm}
&&\qquad\qquad\qquad\qquad\quad{}\cdot[\bm P_A^{r-u-1}\bm P_B^s(\bm Z-\bm\Pi)\bm P_A^r\bm P_B^v]_{jk}(\bm P_B)_{kl}w(k,l) \nonumber\\
&&\qquad{}+\sum_{u=0}^{s-1}\sum_{v=0}^{r-1}\sum_{i,j,k,l}(\bm \pi\bm P_A^r\bm P_B^u)_i(\bm P_B)_{ij}w(i,j) \nonumber\\
\noalign{\vglue-3mm}
&&\qquad\qquad\qquad\qquad\quad{}\cdot[\bm P_B^{s-u-1}(\bm Z-\bm\Pi)\bm P_A^v]_{jk}(\bm P_A)_{kl}w(k,l) \nonumber\\
&&\qquad{}+\sum_{u=0}^{s-1}\sum_{v=0}^{s-1}\sum_{i,j,k,l}(\bm\pi\bm P_A^r\bm P_B^u)_i(\bm P_B)_{ij}w(i,j) \nonumber\\
\noalign{\vglue-3mm}
&&\qquad\qquad\qquad\qquad\quad{}\cdot[\bm P_B^{s-u-1}(\bm Z-\bm\Pi)\bm P_A^r\bm P_B^v]_{jk}(\bm P_B)_{kl}w(k,l),
\end{eqnarray}
where $\bm Z$ is the fundamental matrix associated with $\bm P:=\bm P_A^r\bm P_B^s$.  We conclude that 
\begin{equation}\label{sigma2-periodic}
\sigma^2={1\over r+s}\bigg({\Var}(\eta_1)+2\sum_{m=1}^\infty\Cov(\eta_1,\eta_{m+1})\bigg),
\end{equation}
which relies on (\ref{var}) and (\ref{covsum}).

We summarize the results of this section in the following theorem.

\begin{theorem}
Let $\mu$ and $\sigma^2$ be as in (\ref{mu-periodic}) and (\ref{sigma2-periodic}).
Under the assumptions of the second paragraph of this section, but with the distribution of $X_0$ arbitrary,
$$
\lim_{n\to\infty}n^{-1}\E[S_n]=\mu\quad\text{and}\quad{S_n\over n}\to \mu\;\;{\rm a.s.}
$$
and, if $\sigma^2>0$, 
$$
\lim_{n\to\infty}n^{-1}\Var(S_n)=\sigma^2\quad\text{and}\quad{S_n-n\mu\over\sqrt{n\sigma^2}}\to_d N(0,1).
$$
If $\mu=0$ and $\sigma^2>0$, then (\ref{liminf-limsup}) holds.
\end{theorem}

\begin{remark}
It follows that a pattern is losing, winning, or fair according to whether $\mu<0$, $\mu>0$, or $\mu=0$.
\end{remark}

\begin{proof}
The proof is similar to that of Theorem 1, except that $N:=\min\{n\ge0: X_n=i_0,\;n{\rm\ is\ divisible\ by\ }r+s\}$.
\end{proof}

In the examples to which we will be applying the above mean and variance formulas, additional simplifications will occur because $\bm P_A$ has a particularly simple form and $\bm P_B$ has a spectral representation.  Denote by $\bm P_A'$ the matrix with $(i,j)$th entry $(\bm P_A)_{ij}w(i,j)$, and assume that $\bm P_A'$ \textit{has row sums equal to} 0.  Denote by $\bm P_B'$ the matrix with $(i,j)$th entry $(\bm P_B)_{ij}w(i,j)$, and define $\bm \zeta:=\bm P_B'(1,1,\ldots,1)^\T$ to be the vector of row sums of $\bm P_B'$. Further, let $t:=|\Sigma|$ and suppose that $\bm P_B$ has eigenvalues $1,e_1,\ldots,e_{t-1}$ and corresponding linearly independent right eigenvectors $\bm r_0,\bm r_1,\ldots,\bm r_{t-1}$.  Put $\bm D:={\rm diag}(1,e_1,\ldots,e_{t-1})$ and $\bm R:=(\bm r_0,\bm r_1,\ldots,\bm r_{t-1})$.  Then the rows of $\bm L:=\bm R^{-1}$ are left eigenvectors and $\bm P_B^v=\bm R\bm D^v\bm L$ for all $v\ge0$.  Finally, with
$$
\bm D_v:={\rm diag}\left(v,{1-e_1^v\over1-e_1},\ldots,{1-e_{t-1}^v\over1-e_{t-1}}\right),
$$
we have $\bm I+\bm P_B+\cdots+\bm P_B^{v-1}=\bm R\bm D_v\bm L$ for all $v\ge1$.  Additional notation includes $\bm \pi_{s,r}$ for the unique stationary distribution of $\bm P_B^s\bm P_A^r$, and $\bm\pi$ as before for the unique stationary distribution of $\bm P:=\bm P_A^r\bm P_B^s$.  Notice that $\bm\pi\bm P_A^r=\bm\pi_{s,r}$.  Finally, we assume as well that $w(i,j)=\pm1$ whenever $(\bm P_A)_{ij}>0$ or $(\bm P_B)_{ij}>0$.

These assumptions allow us to write (\ref{mu-periodic}) as
\begin{equation}\label{mu-periodic-simple}
\mu=(r+s)^{-1} {\bm \pi}_{s,r} {\bm R} {\bm D}_s {\bm L}{\bm\zeta},
\end{equation}
and similarly (\ref{var}) and (\ref{covsum}) become
\begin{eqnarray}\label{var-sec5}
\Var(\eta_1)
&=&r+s-\sum_{v=0}^{s-1}({\bm\pi}_{s,r}\bm P_B^v\bm\zeta)^2+2\sum_{u=0}^{r-1}\bm\pi\bm P_A^{u}\bm P^\prime_{A}\bm P^{r-u-1}_A\bm R\bm D_s\bm L\bm\zeta\\
&&\;{}+2\sum_{v=1}^{s-1}\sum_{u=0}^{v-1}\bm\pi_{s,r}\bm P_B^u\bm P^\prime_{B}\bm P_B^{v-u-1}\bm\zeta
-2\sum_{v=1}^{s-1}(\bm\pi_{s,r}\bm R\bm D_v\bm L\bm\zeta)(\bm\pi_{s,r}\bm P_B^v\bm\zeta),\nonumber
\end{eqnarray}
and 
\begin{eqnarray}\label{covsum-sec5}
\sum_{m=1}^\infty\Cov(\eta_1,\eta_{m+1})
&=&\sum_{u=0}^{r-1}\bm\pi\bm P_A^u\bm P^\prime_{A}\bm P_A^{r-u-1}\bm P_B^s(\bm Z-\bm\Pi)\bm P_A^r\bm R\bm D_s\bm L\bm\zeta\\
&&\quad{}+\sum_{u=0}^{s-1}\bm\pi_{s,r}\bm P_B^u\bm P^\prime_{B}\bm P_B^{s-u-1}(\bm Z-\bm\Pi)\bm P_A^r\bm R\bm D_s\bm L\bm\zeta.\nonumber
\end{eqnarray}

A referee has suggested an alternative approach to the results of this section that has certain advantages.  Instead of starting with a time-inhomogeneous Markov chain in $\Sigma$, we begin with the (time-homogeneous) Markov chain in the product space $\bar\Sigma:=\{0,1,\ldots,r+s-1\}\times\Sigma$ with transition probabilities
$$
(i,j)\mapsto(i+1\;{\rm mod\ }r+s,\,k)\textrm{\ \ with\ probability\ }\begin{cases}(\bm P_A)_{jk}&\text{if $0\le i\le r-1$},\\
(\bm P_B)_{jk}&\text{if $r\le i\le r+s-1$}.\end{cases}
$$
It is irreducible and periodic with period $r+s$.
Ordering the states of $\bar\Sigma$ lexicographically, we can write the one-step transition matrix in block form as
$$
\setlength{\arraycolsep}{3mm}
\bar{\bm P}:=\left(\begin{array}{cccccc}
\bm 0 & \bm P_0 & \bm 0 & \cdots & \bm 0\\
\bm 0 & \bm 0 & \bm P_1 & & \bm 0\\
\vdots & \vdots & & \ddots &\\
\bm 0 & \bm 0 & \bm 0 & & \bm P_{r+s-2}\\
\bm P_{r+s-1} & \bm 0 & \bm 0 & \cdots & \bm 0 
\end{array}\right),
$$
where $\bm P_i:=\bm P_A$ for $0\le i\le r-1$ and $\bm P_i:=\bm P_B$ for $r\le i\le r+s-1$.  With $\bm\pi$ as above, the unique stationary distribution for $\bar{\bm P}$ is $\bar{\bm\pi}:=(r+s)^{-1}(\bm\pi_0,\bm\pi_1,\ldots,\bm\pi_{r+s-1})$, where $\bm\pi_0:=\bm\pi$ and $\bm\pi_i:=\bm\pi_{i-1}\bm P_{i-1}$ for $1\le i\le r+s-1$.  Using the idea in (\ref{boldXs}) and (\ref{etas}), we can deduce the strong mixing property and  the central limit theorem.  The mean and variance parameters are as in (\ref{mu,sigma2}) but with bars on the matrices.  Of course, $\bar{\bm Z}:=(\bm I-(\bar{\bm P}-\bar{\bm\Pi}))^{-1}$, $\bar{\bm\Pi}$ being the square matrix each of whose rows is $\bar{\bm\pi}$; $\bar{\bm P}'$ is as $\bar{\bm P}$ but with $\bm P_A'$ and $\bm P_B'$ in place of $\bm P_A$ and $\bm P_B$; and similarly for $\bar{\bm P}''$.

The principal advantage of this approach is the simplicity of the formulas for the mean and variance parameters.  Another advantage is that patterns other than those of the form $[r,s]$ can be treated just as easily.  The only drawback is that, in the context of our two models, the matrices are no longer $3\times3$ or $4\times4$ but rather $3(r+s)\times3(r+s)$ and $4(r+s)\times4(r+s)$.  
In other words, the formulas (\ref{mu-periodic-simple})--(\ref{covsum-sec5}), although lacking elegance, are more user-friendly.

\section{Patterns of capital-dependent games}

Let games $A$ and $B$ be as in Section 3; see especially (\ref{probs-profit}).  Both games are losing.  In this section we show that, for every $\rho\in(0,1)$, and for every pair of positive integers $r$ and $s$ except for $r=s=1$, the pattern $[r,s]$, which stands for $r$ plays of game $A$ followed by $s$ plays of game $B$, is winning for sufficiently small $\eps>0$.  Notice that it will suffice to treat the case $\eps=0$.

We begin by finding a formula for $\mu_{[r,s]}(0)$, the asymptotic mean per game played of the player's cumulative profit for the pattern $[r,s]$, assuming $\eps=0$.  Recall that $\bm P_A$ is given by (\ref{P-profit}) with $p_0=p_1=p_2:={1\over2}$, and $\bm P_B$ is given by (\ref{P-profit}) with $p_0:=\rho^2/(1+\rho^2)$ and $p_1=p_2:=1/(1+\rho)$, where $\rho>0$. First, we can prove by induction that
$$
\setlength{\arraycolsep}{1.5mm}
{\bm P}_A^r=\left(\begin{array}{ccc}
1-2a_r&a_r&a_r\\
a_r&1-2a_r&a_r\\
a_r&a_r&1-2a_r
\end{array}\right),\qquad r\ge1,
$$
where 
\begin{equation}\label{a_r}
a_r:={1-(-{1\over2})^r\over3}.
\end{equation}
Next, with $S:=\sqrt{(1+\rho^2)(1+4\rho+\rho^2)}$, the nonunit eigenvalues of ${\bm P}_B$ are
$$
e_1:=-{1\over2}+{(1-\rho)S\over2(1+\rho)(1+\rho^2)},\qquad
e_2:=-{1\over2}-{(1-\rho)S\over2(1+\rho)(1+\rho^2)},
$$
and we define the diagonal matrix $\bm D:={\rm diag}(1,e_1,e_2)$. 
The corresponding right eigenvectors
\begin{eqnarray*}
{\bm r}_0:=\left(\begin{array}{c}
1\\
1\\
1
\end{array}\right),\quad
{\bm r}_1&:=&\left(\begin{array}{c}
(1+\rho)(1-\rho^2-S)\\
2+\rho+2\rho^2+\rho^3+\rho S\\
-(1+2\rho+\rho^2+2\rho^3-S)
\end{array}\right),\\
{\bm r}_2&:=&\left(\begin{array}{c}
(1+\rho)(1-\rho^2+S)\\
2+\rho+2\rho^2+\rho^3-\rho S\\
-(1+2\rho+\rho^2+2\rho^3+S)
\end{array}\right),
\end{eqnarray*}
are linearly independent for all $\rho>0$ (including $\rho=1$), 
so we define ${\bm R}:=({\bm r}_0,{\bm r}_1,{\bm r}_2)$ and ${\bm L}:={\bm R}^{-1}$.  

Then
$$
{\bm P}_B^s={\bm R}{\bm D}^s{\bm L},\qquad s\ge0,
$$
which leads to an explicit formula for ${\bm P}_B^s{\bm P}_A^r$, from which we can compute its unique stationary distribution ${\bm \pi}_{s,r}$ as a left eigenvector corresponding to the eigenvalue 1. With 
$$
{\bm\zeta}:=(2\rho^2/(1+\rho^2)-1,\,2/(1+\rho)-1,\,2/(1+\rho)-1)^\T,
$$
we can use (\ref{mu-periodic-simple}) to write
$$
\mu_{[r,s]}(0)={1\over r+s}\,{\bm \pi}_{s,r}{\bm R}\,\,{\rm diag}\!\left(s,\,{1-e_1^s\over1-e_1},\,{1-e_2^s\over1-e_2}\right){\bm L}{\bm\zeta}.
$$
Algebraic computations show that this reduces to 
\begin{equation}\label{E/D}
\mu_{[r,s]}(0)=E_{r,s}/D_{r,s},
\end{equation}
where
\begin{eqnarray*}
E_{r,s}&:=&3a_r\{[2+(3a_r-1)(e_1^s+e_2^s-2e_1^s e_2^s)-(e_1^s+e_2^s)](1-\rho)(1+\rho)S\\
&&\qquad\qquad{}+a_r(e_2^s-e_1^s)[5(1+\rho)^2(1+\rho^2)-4\rho^2]\}(1-\rho)^2
\end{eqnarray*}
and
$$
D_{r,s}:=4(r+s)[1+(3a_r-1)e_1^s][1+(3a_r-1)e_2^s](1+\rho+\rho^2)^2S.
$$
This formula will allow us to determine the sign of $\mu_{[r,s]}(0)$ for all $r$, $s$, and $\rho$.

\begin{theorem}
Let games $A$ and $B$ be as in Theorem 2 (with the bias parameter absent).  For every pair of positive integers $r$ and $s$ except $r=s=1$, $\mu_{[r,s]}(0)>0$ for all $\rho\in(0,1)$, $\mu_{[r,s]}(0)=0$ for $\rho=1$, and $\mu_{[r,s]}(0)<0$ for all $\rho>1$.  In addition, $\mu_{[1,1]}(0)=0$ for all $\rho>0$.
\end{theorem}

The last assertion of the theorem was known to Pyke (2003).  As before, the Parrondo effect appears (with the bias parameter present) if and only if $\mu_{[r,s]}(0)>0$.  A reverse Parrondo effect, in which two winning games combine to lose, appears (with a \textit{negative} bias parameter present) if and only if $\mu_{[r,s]}(0)<0$.

\begin{corollary}
Let games $A$ and $B$ be as in Corollary 3 (with the bias parameter present).
For every $\rho\in(0,1)$, and for every pair of positive integers $r$ and $s$ except $r=s=1$, there exists $\eps_0>0$, depending on $\rho$, $r$, and $s$, such that Parrondo's paradox holds for games $A$ and $B$ and pattern $[r,s]$, that is, $\mu_A(\eps)<0$, $\mu_B(\eps)<0$, and $\mu_{[r,s]}(\eps)>0$, whenever $0<\eps<\eps_0$. 
\end{corollary}

\begin{proof}[Proof of theorem]
Temporarily denote $\mu_{[r,s]}(0)$ by $\mu_{[r,s]}(\rho,0)$ to emphasize its dependence on $\rho$.  Then, given $\rho\in(0,1)$, the same argument that led to (\ref{symmetry}) yields
\begin{equation}\label{mu=-mu}
\mu_{[r,s]}(1/\rho,0)=-\mu_{[r,s]}(\rho,0), \qquad r,s\ge1.
\end{equation}
(This also follows from (\ref{E/D}).)
It is therefore enough to treat the case $\rho\in(0,1)$.  Notice that $3a_r-1=-(-{1\over2})^r$ for all $r\ge1$ 
and $e_1,e_2\in(-1,0)$, hence $D_{r,s}>0$ for all integers $r,s\ge1$.  To show that $\mu_{[r,s]}(0)>0$ it will suffice to show that $E_{r,s}>0$, which is equivalent to showing that
\begin{eqnarray}\label{Efactor}
&&2+(3a_r-1)(e_1^s+e_2^s-2e_1^s e_2^s)-(e_1^s+e_2^s)\nonumber\\
&&\qquad\qquad\quad{}+a_r(e_2^s-e_1^s){5(1+\rho)^2(1+\rho^2)-4\rho^2\over(1-\rho)(1+\rho)S}
\end{eqnarray}
is positive. We will show that, except when $r=s=1$, (\ref{Efactor}) is positive for all integers $r,s\ge1$.

First we consider the case in which $r\ge2$ and $s=1$.  Noting that $e_1+e_2=-1$, $e_2-e_1=-(1-\rho)S/[(1+\rho)(1+\rho^2)]$, and $e_1e_2=2\rho^2/[(1+\rho)^2(1+\rho^2)]$, (\ref{Efactor}) becomes
\begin{eqnarray*}
&&3+\bigg(\!\!-{1\over2}\bigg)^r\bigg(1+{4\rho^2\over(1+\rho)^2(1+\rho^2)}\bigg)-{1-(-{1\over2})^r\over3}\,{5(1+\rho)^2(1+\rho^2)-4\rho^2\over(1+\rho)^2(1+\rho^2)}\\
&&\quad{}>3-{1\over8}\bigg(1+{4\rho^2\over(1+\rho)^2(1+\rho^2)}\bigg)-{3\over8}\bigg(5-{4\rho^2\over(1+\rho)^2(1+\rho^2)}\bigg)\nonumber\\
&&\quad{}=1+{\rho^2\over(1+\rho)^2(1+\rho^2)}>0.
\end{eqnarray*}
On the other hand, if $r=s=1$, the left side of the last equation becomes
$$
3-{1\over2}\bigg(1+{4\rho^2\over(1+\rho)^2(1+\rho^2)}\bigg)-{1\over2}\bigg(5-{4\rho^2\over(1+\rho)^2(1+\rho^2)}\bigg)=0,
$$
so $\mu_{[1,1]}(0)=0$.   

Next we consider the case of $r\ge2$ and $s$ odd, $s\ge3$.  Since $\rho\in(0,1)$, we have $e_1=-{1\over2}(1-x)$ and $e_2=-{1\over2}(1+x)$, where $0<x<1$.  Therefore, since $s$ is odd,
\begin{eqnarray}\label{e2-e1}
e_2^s-e_1^s&=&-[(1+x)^s-(1-x)^s]/2^s\nonumber\\
&=&-{1\over2^s}\sum_{k=0}^s{s\choose k}[1-(-1)^k]x^k \nonumber\\
&=&-x\,{1\over2^{s-1}}\sum_{1\le k\le s;\; k{\rm\ odd}}{s\choose k}x^{k-1} \nonumber\\
&>&(e_2-e_1){1\over2^{s-1}}\sum_{1\le k\le s;\; k{\rm\ odd}}{s\choose k} \nonumber\\
&=&e_2-e_1,
\end{eqnarray}
where the last identity uses the fact that, when $s$ is odd, the odd-numbered binomial coefficients ${s\choose1},{s\choose3},\ldots,{s\choose s}$ are equal to the even-numbered ones ${s\choose s-1},\break{s\choose s-3},\ldots,{s\choose 0}$, hence the sum of the odd-numbered ones is $2^{s-1}$.

It follows that, since the case $r=1$ is excluded,
\begin{eqnarray*}
&&2+(3a_r-1)(e_1^s+e_2^s-2e_1^s e_2^s)-(e_1^s+e_2^s)\\
&&\quad\qquad{}+a_r(e_2^s-e_1^s){5(1+\rho)^2(1+\rho^2)-4\rho^2\over(1-\rho)(1+\rho)S}\\
&&\quad{}>2+{1\over8}(e_1+e_2-2e_1 e_2)+{3\over8}(e_2-e_1){5(1+\rho)^2(1+\rho^2)-4\rho^2\over(1-\rho)(1+\rho)S}\\
&&\quad{}=2-{1\over8}\bigg(1+{4\rho^2\over(1+\rho)^2(1+\rho^2)}\bigg)-{3\over8}\bigg(5-{4\rho^2\over(1+\rho)^2(1+\rho^2)}\bigg)\\
&&\quad{}={\rho^2\over(1+\rho)^2(1+\rho^2)}>0.
\end{eqnarray*}

Consider next the case in which $r\ge1$ and $s$ is even.  In this case, the last term in (\ref{Efactor}) is positive, and the remaining terms are greater than $2-(1/4)(5/4)-5/4>0$.

Next we treat the case $[1,3]$ separately.  Our formula (\ref{E/D}) gives
$\mu_{[1,3]}(0)=3\rho^2(1-\rho)^3(1+\rho)(1+2\rho+2\rho^3+\rho^4)/[4(1+6\rho+24\rho^2+62\rho^3+111\rho^4+156\rho^5+180\rho^6+156\rho^7+111\rho^8+62\rho^9+24\rho^{10}+6\rho^{11}+\rho^{12})]$, which is positive since $\rho\in(0,1)$.

It remains to consider the case of $r=1$ and $s$ odd, $s\ge5$.  We must show that
$$
2-{1\over2}(e_1^s+e_2^s+2e_1^s e_2^s)+{1\over2}(e_2^s-e_1^s){5(1+\rho)^2(1+\rho^2)-4\rho^2\over(1-\rho)(1+\rho)S}>0.
$$
Here the inequality of (\ref{e2-e1}) is too crude.  Again, we have $e_1=-{1\over2}(1-x)$ and $e_2=-{1\over2}(1+x)$, where $x=(1-\rho)S/[(1+\rho)(1+\rho^2)]\in(0,1)$, so we can replace $(1-\rho)S$ by $x(1+\rho)(1+\rho^2)$.  We can eliminate the $-4\rho^2$ term in the fraction, and we can replace the $2e_1^se_2^s$ term by $-2e_1^s$; in the first case we are eliminating a positive contribution, and in the second case we are making a negative contribution more negative.  Thus, it will suffice to prove that
$$
2^{s+2}+[(1+x)^s-(1-x)^s]-[(1+x)^s-(1-x)^s](5/x)>0
$$
for $0<x<1$, or equivalently that
$$
f(x):={2^{s+2}x\over(1+x)^s-(1-x)^s}>5-x
$$
for $0<x<1$. Now $[(1+x)^s-(1-x)^s]/x$ is a polynomial of degree $s-1$ with nonnegative coefficients, so $f$ is decreasing on $(0,1)$ with $f(0+)=2^{s+1}/s$ and $f(1-)=4$.  Further, $f({1\over2})=2^{s+1}/[(3/2)^s-(1/2)^s]>2^{s+1}/(3/2)^s=2(4/3)^s>5$ for $s\ge5$.  Thus, our inequality $f(x)>5-x$ holds for $0<x\le{1\over2}$.  It remains to confirm it for ${1\over2}<x<1$.  On that interval we claim that
$f(x)+x$ is decreasing, and it clearly approaches 5 as $x\to1-$.  Thus, it remains to show that $f'(x)<-1$ for ${1\over2}<x<1$.  We compute
\begin{eqnarray*}
f'(x)&=&2^{s+2}{[(1+x)^s-(1-x)^s]-sx[(1+x)^{s-1}+(1-x)^{s-1}]\over[(1+x)^s-(1-x)^s]^2}\\
&<&2^{s+2}{(1+x)^s-sx(1+x)^{s-1}\over[(1+x)^s-(1-x)^s]^2}\\
&=&-2^{s+2}{[(s-1)x-1](1+x)^{s-1}\over[(1+x)^s-(1-x)^s]^2}\\
&<&-4{[4(1/2)-1](1+x)^{s-1}\over(1+x)^s-(1-x)^s}\\
&<&-2{(1+x)^s\over(1+x)^s-(1-x)^s}\\
&<&-2,
\end{eqnarray*}
as required.
\end{proof}

Here are four of the simplest cases:
\begin{eqnarray*}
\mu_{[1,1]}(\eps)&=&-{3(1+2\rho+18\rho^2+2\rho^3+\rho^4)\over4(1+\rho+\rho^2)^2}\,\eps+O(\eps^2),\\
\mu_{[1,2]}(\eps)&=&{(1-\rho)^3(1+\rho)(1+2\rho+\rho^2+2\rho^3+\rho^4)\over3+12\rho+20\rho^2+28\rho^3+36\rho^4+28\rho^5+20\rho^6+12\rho^7+3\rho^8}+O(\eps),\\
\mu_{[2,1]}(\eps)&=&{(1-\rho)^3(1+\rho)\over10+20\rho+21\rho^2+20\rho^3+10\rho^4}+O(\eps),\\
\mu_{[2,2]}(\eps)&=&{3(1-\rho)^3(1+\rho)\over8(3+6\rho+7\rho^2+6\rho^3+3\rho^4)}+O(\eps).
\end{eqnarray*}

We turn finally to the evaluation of the asymptotic variance per game played of the player's cumulative profit.  We denote this variance by $\sigma^2_{[r,s]}(\eps)$, and we note that it suffices to assume that $\eps=0$ to obtain the formula up to $O(\eps)$.

A formula for $\sigma_{[r,s]}^2(0)$ analogous to (\ref{E/D}) would be extremely complicated.  We therefore consider the matrix formulas of Section 5 to be in final form.

Temporarily denote $\sigma_{[r,s]}^2(0)$ by $\sigma_{[r,s]}^2(\rho,0)$ to emphasize its dependence on $\rho$.  Then, given $\rho\in(0,1)$, the argument that led to (\ref{mu=-mu}) yields
$$
\sigma_{[r,s]}^2(1/\rho,0)=\sigma_{[r,s]}^2(\rho,0), \qquad r,s\ge1.
$$

For example, we have
\begin{eqnarray*}
\sigma^2_{[1,1]}(\eps)&=&\left({3\rho}\over{1+\rho+\rho^2}\right)^2+O(\eps),\\
\sigma^2_{[1,2]}(\eps)&=&3(8+104\rho+606\rho^2+2220\rho^3+6189\rho^4+14524\rho^5+29390\rho^6\\
&&\quad\;{}+51904\rho^7+81698\rho^8+115568\rho^9+147156\rho^{10}+169968\rho^{11}\\
&&\quad\;{}+178506\rho^{12}+169968\rho^{13}+147156\rho^{14}+115568\rho^{15}+81698\rho^{16}\\
&&\quad\;{}+51904\rho^{17}+29390\rho^{18}+14524\rho^{19}+6189\rho^{20}+2220\rho^{21}\\
&&\quad\;{}+606\rho^{22}+104\rho^{23}+8\rho^{24})\\
&&\;\;/(3+12\rho+20\rho^2+28\rho^3+36\rho^4+28\rho^5+20\rho^6+12\rho^7+3\rho^8)^3\\
&&\quad\;{}+O(\eps),\\
\sigma^2_{[2,1]}(\eps)&=&3(414+2372\rho+6757\rho^2+13584\rho^3+21750\rho^4+28332\rho^5+30729\rho^6\\
&&\quad\;{}+28332\rho^7+21750\rho^8+13584\rho^9+6757\rho^{10}+2372\rho^{11}+414\rho^{12})\\
&&\;\;/(10+20\rho+21\rho^2+20\rho^3+10\rho^4)^3+O(\eps),\\
\sigma^2_{[2,2]}(\eps)&=&9(25+288\rho+1396\rho^2+4400\rho^3+10385\rho^4+19452\rho^5+29860\rho^6\\
&&\quad\;{}+38364\rho^7+41660\rho^8+38364\rho^9+29860\rho^{10}+19452\rho^{11}\\
&&\quad\;{}+10385\rho^{12}+4400\rho^{13}+1396\rho^{14}+288\rho^{15}+25\rho^{16})\\
&&\;\;/[16(1+\rho+\rho^2)^2(3+6\rho+7\rho^2+6\rho^3+3\rho^4)^3]+O(\eps).
\end{eqnarray*}
As functions of $\rho\in(0,1)$, $\sigma_{[1,1]}^2(0)$, $\sigma_{[1,2]}^2(0)$, and $\sigma_{[2,2]}^2(0)$ are increasing, whereas $\sigma_{[2,1]}^2(0)$ is decreasing.  All approach 1 as $\rho\to1-$.  Their limits as $\rho\to0+$ are respectively 0, 8/9, 25/48, and 621/500.

It should be explicitly noted that 
\begin{equation}\label{coincidence}
\mu_{[1,1]}(0)=\mu_B(0)\quad{\rm and}\quad\sigma_{[1,1]}^2(0)=\sigma_B^2(0)
\end{equation}
for all $\rho>0$.  Since these identities may be counterintuitive, let us derive them in greater detail.  We temporarily denote the stationary distributions and fundamental matrices of $\bm P_B$, $\bm P_A\bm P_B$, and $\bm P_B\bm P_A$, with subscripts $B$, $AB$, and $BA$, respectively.  Then
$$
\bm\pi_B={1\over2(1+\rho+\rho^2)}(1+\rho^2,\,\rho(1+\rho),\,1+\rho)
$$
and
$$
\bm\pi_{BA}={1\over2(1+\rho+\rho^2)}(1+\rho^2,\,1+\rho,\,\rho(1+\rho)),
$$
so
$$
\mu_B(0)=\bm\pi_B\bm\zeta=0\quad{\rm and}\quad\mu_{[1,1]}(0)={1\over2}\bm\pi_{BA}\bm\zeta=0.
$$
As for the variances, we recall that
$$
\sigma^2_B(0)=1+2\bm\pi_B\bm P_B^\prime\bm Z_B\bm\zeta=1+\bm\pi_B\bm P_B^\prime(2\bm Z_B\bm\zeta)
$$
and
$$
\sigma^2_{[1,1]}(0)=1+\bm\pi_{AB}\bm P_A^\prime(\bm I+\bm P_B\bm Z_{AB}\bm P_A)\bm\zeta.
$$
But $\bm\pi_{AB}\bm P_A^\prime=\bm\pi_{BA}\bm P_B\bm P_A^\prime$ and
\begin{eqnarray*}
\bm\pi_B\bm P_B^\prime &=&{1-\rho\over2(1+\rho+\rho^2)}(1+\rho,\,-\rho,\,-1),
\\
\bm\pi_{BA}\bm P_B\bm P_A^\prime&=& 
{1-\rho\over2(1+\rho+\rho^2)}(1+\rho,\,-1,\,-\rho),
\\
2\bm Z_B\bm\zeta&=&{1-\rho\over(1+\rho+\rho^2)^2}\left(\begin{array}{c}
-(1+\rho)(1+4\rho+\rho^2)\\
3+\rho+\rho^2+\rho^3\\
1+\rho+\rho^2+3\rho^3
\end{array}\right),\\
(\bm I +\bm P_B\bm Z_{AB}\bm P_A)\bm\zeta&=&{1-\rho\over(1+\rho+\rho^2)^2}
\left(\begin{array}{c}
-(1+\rho)(1+4\rho+\rho^2)\\
1+\rho+\rho^2+3\rho^3\\
3+\rho+\rho^2+\rho^3
\end{array}\right),
\end{eqnarray*}
so it follows that
$$
\sigma^2_B(0)=\sigma^2_{[1,1]}(0)=\bigg({3\rho\over1+\rho+\rho^2}\bigg)^2.
$$

An interesting special case of (\ref{coincidence}) is the case $\rho=0$.  Technically, we have ruled out this case because we want our underlying Markov chain to be irreducible.  Nevertheless, the games are well defined.  Assuming $X_0=0$, repeated play of game $B$ leads to the deterministic sequence
$$
(S_1,S_2,S_3,\ldots)=(-1,0,-1,0,-1,0,\ldots),
$$
hence $\mu_B(0)=\lim_{n\to\infty}n^{-1}\E[S_n]=0$ and $\sigma_B^2(0)=\lim_{n\to\infty}n^{-1}\Var(S_n)=0$.  On the other hand, repeated play of pattern $AB$ leads to the random sequence
$$
(S_1,S_2,S_3,\ldots)=(-1,0,-1,0,\ldots,-1,0,1,2,2\pm1,2,2\pm1,2,2\pm1,\ldots),
$$
where the number $N$ of initial $(-1,0)$ pairs is the number of losses at game $A$ before the first win at game $A$ (in particular, $N$ is nonnegative geometric$({1\over2})$), and the $\pm1$ terms signify independent random variables that are $\pm1$ with probability ${1\over2}$ each and independent of $N$.  Despite the randomness, the sequence is bounded, so we still have $\mu_{[1,1]}(0)=0$ and $\sigma_{[1,1]}^2(0)=0$.

When $\rho = 1/3$, the mean and variance formulas simplify to
\begin{align*}
\mu_{[1,1]}(\eps)&=-{228\over169}\,\eps+O(\eps^2),&\sigma^2_{[1,1]}(\eps)&={81\over169}+O(\eps),\\
\mu_{[1,2]}(\eps)&={2416\over35601}+O(\eps),&\sigma^2_{[1,2]}(\eps)&={14640669052339\over15040606062267}+O(\eps),\\
\mu_{[2,1]}(\eps)&={32\over1609}+O(\eps),&\sigma^2_{[2,1]}(\eps)&={4628172105\over4165509529} +O(\eps),\\
\mu_{[2,2]}(\eps)&={4\over163}+O(\eps),&\sigma^2_{[2,2]}(\eps)&={1923037543\over2195688729}+O(\eps).
\end{align*}
Pyke (2003) obtained $\mu_{[2,2]}(0)\approx0.0218363$ when $\rho=1/3$, but that number is inconsistent with Ekhad and Zeilberger's (2000) calculation, confirmed above,  that $\mu_{[2,2]}(0)=4/163\approx0.0245399$.

\section{Patterns of history-dependent games}

Let games $A$ and $B$ be as in Section 4; see especially (\ref{probs-history}).  Both games are losing.  In this section we attempt to find conditions on $\kappa$ and $\lambda$ such that, for every pair of positive integers $r$ and $s$, the pattern $[r,s]$, which stands for $r$ plays of game $A$ followed by $s$ plays of game $B$, is winning for sufficiently small $\eps>0$.  Notice that it will suffice to treat the case $\eps=0$.

We begin by finding a formula for $\mu_{[r,s]}(0)$, the asymptotic mean per game played of the player's cumulative profit for the pattern $[r,s]$, assuming $\eps=0$.  Recall that $\bm P_A$ is defined by (\ref{P-history}) with $p_0=p_1=p_2=p_3:={1\over2}$, and $\bm P_B$ is defined by (\ref{P-history}) with $p_0:=1/(1+\kappa)$, $p_1=p_2:=\lambda/(1+\lambda)$, and $p_3:=1-\lambda/(1+\kappa)$, where $\kappa>0$, $\lambda>0$, and $\lambda<1+\kappa$.  First, it is clear that ${\bm P}_A^r={\bm U}$ for all $r\ge2$, where $\bm U$ is the $4\times4$ matrix with all entries equal to 1/4.  

As for the spectral representation of $\bm P_B$, its eigenvalues include 1 and the three roots of the cubic equation
$x^3+a_2x^2+a_1x+a_0=0$, where
\begin{equation}\label{cubic-coeffs}
a_2:={\lambda-\kappa\over1+\kappa},\quad
a_1:={(\lambda-\kappa)\lambda(2+\kappa+\lambda)\over(1+\kappa)^2(1+\lambda)^2},\quad
a_0:=-{(1-\kappa\lambda)(1+\kappa-\lambda-\lambda^2)\over(1+\kappa)^2(1+\lambda)^2}.
\end{equation}
With the help of Cardano's formula, we find that
the nonunit eigenvalues of ${\bm P}_B$ are
\begin{eqnarray*}
e_1&:=&\frac{P+Q}{3(1+\kappa)(1+\lambda)}-\frac{\lambda-\kappa}{3(1+\kappa)},\\
e_2&:=&\frac{\omega P+\omega^2 Q}{3(1+\kappa)(1+\lambda)}-\frac{\lambda-\kappa}{3(1+\kappa)},\\
e_3&:=&\frac{\omega^2 P+\omega Q}{3(1+\kappa)(1+\lambda)}-\frac{\lambda-\kappa}{3(1+\kappa)}, 
\end{eqnarray*}
where $\omega:=e^{2\pi i/3}=-{1\over2}+{1\over2}\sqrt{3}\, i$ and $\omega^2=e^{4\pi i/3}=\bar\omega=-{1\over2}-{1\over2}\sqrt{3}\, i$ are cube roots of unity, and 
$$
P:=\sqrt[3]{{\beta+\sqrt{\beta^2+ 4\alpha^3}\over2}},\qquad
Q:=\sqrt[3]{{\beta-\sqrt{\beta^2+ 4\alpha^3}\over2}},
$$
with 
\begin{eqnarray*}
\alpha&:=&(\lambda-\kappa)(\kappa+5\lambda+5\kappa\lambda+\lambda^2+\kappa\lambda^2-\lambda^3),\\
\beta&:=&(1+\lambda)(27+54\kappa-27\lambda+27\kappa^2-54\kappa\lambda-27\lambda^2+2\kappa^3-42\kappa^2\lambda\\
&&\qquad\qquad{}-30\kappa\lambda^2+16\lambda^3-14\kappa^3\lambda+6\kappa^2\lambda^2+30\kappa\lambda^3+5\lambda^4\\
&&\qquad\qquad{}+2\kappa^3\lambda^2+21\kappa^2\lambda^3+6\kappa\lambda^4-2\lambda^5).
\end{eqnarray*}
Of course, $e_1$, $e_2$, and $e_3$ are each less than 1 in absolute value.

Notice that the definitions of $P$ and $Q$ are slightly ambiguous, owing to the nonuniqueness of the cube roots. (If $(P,Q)$ is replaced in the definitions of $e_1$, $e_2$, and $e_3$ by $(\omega P,\omega^2 Q)$ or by $(\omega^2 P,\omega Q)$, then $e_1$, $e_2$, and $e_3$ are merely permuted.)  If $\beta^2+4\alpha^3>0$, then $P$ and $Q$ can be taken to be real and distinct,\footnote{Caution should be exercised when evaluating $Q$ numerically.  Specifically, if $x<0$, \textit{Mathematica} returns a complex root for $\sqrt[3]{x}$.  If the real root is desired, one should enter $-\sqrt[3]{-x}$.  
This issue does not arise with $P$ and can be avoided with $Q$ by redefining $Q:=-\alpha/P$.} in which case $e_1$ is real and $e_2$ and $e_3$ are complex conjugates; in particular, $e_1$, $e_2$, and $e_3$ are distinct.  If $\beta^2+4\alpha^3=0$, then $P$ and $Q$ can be taken to be real and equal, in which case $e_1$, $e_2$, and $e_3$ are real with $e_2=e_3$.  If $\beta^2+4\alpha^3<0$, then $P$ and $Q$ can be taken to be complex conjugates, in which case $e_1$, $e_2$, and $e_3$ are real and distinct;  in fact, they can be written
\begin{eqnarray*} 
e_1&:=&\frac{2(\sqrt{-\alpha}\,)\cos(\theta/3)}{3(1+\kappa)(1+\lambda)}-\frac{\lambda-\kappa}{3(1+\kappa)},\\
e_2&:=&\frac{2(\sqrt{-\alpha}\,)\cos((\theta+2\pi)/3)}{3(1+\kappa)(1+\lambda)}-\frac{\lambda-\kappa}{3(1+\kappa)},\\
e_3&:=&\frac{2(\sqrt{-\alpha}\,)\cos((\theta+4\pi)/3)}{3(1+\kappa)(1+\lambda)}-\frac{\lambda-\kappa}{3(1+\kappa)}, 
\end{eqnarray*}
where $\theta:=\cos^{-1}({1\over2}\beta/\sqrt{-\alpha^3}\,)$, which implies that $1>e_1>e_3>e_2>-1$.  See Figure 1.

\begin{figure}[htb]
\centering
\includegraphics[width = 250bp]{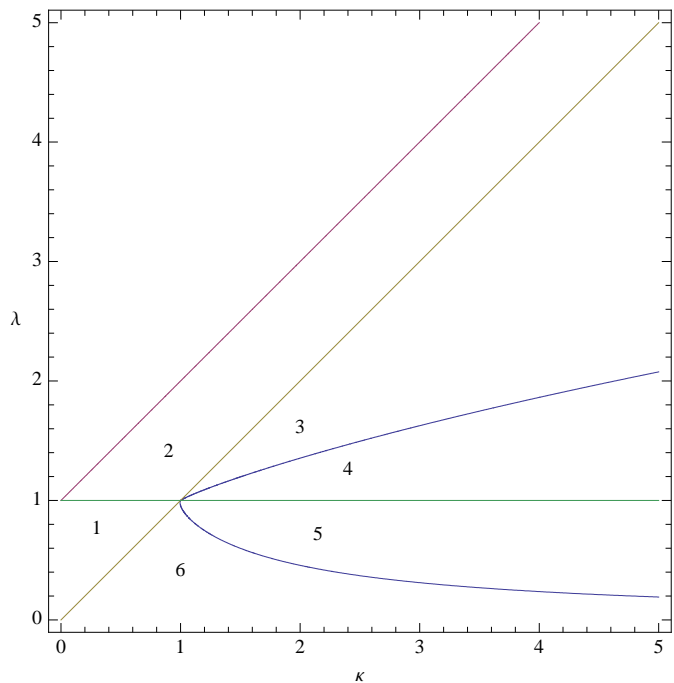}
\caption{The parameter space $\{(\kappa,\lambda): \kappa>0,\;\lambda>0,\;\lambda<1+\kappa\}$, restricted to $\kappa<5$ and $\lambda<5$.  In regions 1, 2, 3, and 6, $\beta^2+4\alpha^3>0$ ($e_1$ is real, $e_2$ and $e_3$ are complex conjugates); in regions 4 and 5, $\beta^2+4\alpha^3<0$ ($e_1$, $e_2$, and $e_3$ are real).  If the conjecture stated below is correct, then, in regions 1, 3, and 4, $\mu_{[r,s]}(0)>0$ for all $r,s\ge1$; in regions 2, 5, and 6, $\mu_{[r,s]}(0)<0$  for all $r,s\ge1$.}
\end{figure}  

If we define the vector-valued function $\bm r : {\bf C} \mapsto {\bf C}^4$ by
\begin{eqnarray*}
&&\bm r(x) \\
&&\;{}:=\left(\begin{array}{c}
-\lambda(1+\kappa-\lambda-\lambda^2)-(1+\lambda)(1+\kappa-\lambda-\lambda^2)x+(1+\kappa)(1+\lambda)^2x^2\\
(1+\kappa-\lambda-\lambda^2)-(1+\kappa)(1+\lambda)(\lambda-\kappa)x-(1+\kappa)^2(1+\lambda)x^2\\
-(1+\kappa-\lambda)(1-\kappa\lambda)+(1+\kappa)(1-\kappa\lambda)x\\
\lambda(1-\kappa\lambda)
\end{array}\right),
\end{eqnarray*}
then 
\begin{eqnarray*}
{\bm r}_0:=(1,\,1,\,1,\,1)^\T,\quad{\bm r}_1:=\bm r(e_1),\quad{\bm r}_2:=\bm r(e_2),\quad{\bm r}_3:=\bm r(e_3),
\end{eqnarray*}
are corresponding right eigenvectors, and they are linearly independent if and only if the eigenvalues are distinct.  If the eigenvalues are distinct, we define ${\bm R}:=({\bm r}_0,{\bm r}_1,{\bm r}_2,{\bm r}_3)$ and ${\bm L}:={\bm R}^{-1}$.

With
$$
{\bm u}:=\frac{1}{4}\left(1,\,1,\,1,\,1\right)\quad{\rm and}\quad
{\bm\zeta}:=\left(\begin{array}{c}
2/(1+\kappa)-1\\
2\lambda/(1+\lambda)-1\\
2\lambda/(1+\lambda)-1\\
2[1-\lambda/(1+\kappa)]-1
\end{array}\right),
$$
we can use (\ref{mu-periodic-simple}) to write 
\begin{equation}\label{mu[r,s]}
\mu_{[r,s]}(0)= (r+s)^{-1}E_s, \qquad r\ge2,\quad s\ge1,
\end{equation}
where
$$
E_{s}:={\bm u}{\bm R}\bm D_s{\bm L}{\bm\zeta}
$$
with
$$
\bm D_s:={\rm diag}\!\left(s,\,\frac{1-e_1^s}{1-e_1},\,\frac{1-e_2^s}{1-e_2},\,\frac{1-e_3^s}{1-e_3}\right).
$$
Algebraic simplification leads to
\begin{equation} \label{Es}
E_s=c_0-c_1 e_1^s-c_2 e_2^s-c_3 e_3^s,
\end{equation}
where
$$
c_0:=\frac{(1+\kappa)(\lambda-\kappa)(1-\lambda)}{4\lambda(2+\kappa+\lambda)},
$$
and
$$
c_1:=f(e_1,e_2,e_3),\quad c_2:=f(e_2,e_3,e_1),\quad c_3:=f(e_3, e_1, e_2),
$$
with
\begin{eqnarray*}
f(x,y,z)&:=&(\lambda-\kappa)[\lambda(\lambda-\kappa)-(1+\kappa-\lambda-\lambda^2)x+(1+\kappa)(1+\lambda)x^2]\\
&&\;\;\;{}\cdot[1+3\kappa-2\lambda+3\kappa^2-4\kappa\lambda-\lambda^2+\kappa^3-9\kappa\lambda^2+6\lambda^3+2\kappa^3\lambda\\
&&\qquad{}-7\kappa^2\lambda^2+6\kappa\lambda^3+\kappa^3\lambda^2-2\kappa^2\lambda^3+4\kappa\lambda^4-2\lambda^5\\
&&\qquad{}-(1+\kappa)(1+\lambda)(1+2\kappa-3\lambda+\kappa^2-2\kappa\lambda+\kappa^2\lambda-2\kappa\lambda^2\\
&&\qquad{}+2\lambda^3)(y+z)+(1+\kappa)^2(1+\lambda)^2(1+\kappa-2\lambda)yz]\\
&&\;\;\;{}/[4(1+\kappa)^3\lambda(1+\lambda)^2(1-\kappa\lambda)(1-x)(x-y)(x-z)].
\end{eqnarray*}

It remains to consider
\begin{equation}\label{mu[1,s]}
\mu_{[1,s]}(0)=(1+s)^{-1}F_s,\qquad s\ge1,
\end{equation}
where
\begin{eqnarray}
\nonumber
F_{s} &:=& {\bm \pi}_{s,1} {\bm R} \bm D_s{\bm L}{\bm\zeta}.
\end{eqnarray}
Here ${\bm \pi}_{s,1}$ is the stationary distribution of $\bm P_B^s \bm P_A$. This is just a left eigenvector of
$$
{\bm R}\,{\rm diag}(1,e_1^s,e_2^s,e_3^s){\bm L}\bm P_A
$$
corresponding to eigenvalue 1. Writing $\hat{\bm \pi}=(\hat{\pi}_0,\hat{\pi}_1,\hat{\pi}_2,\hat{\pi}_3):=\bm\pi_{s,1}$, we find that
\begin{eqnarray*}
\hat{\pi}_0&=&\hat{\pi}_1=\frac{1+f_0(e_1,e_2,e_3)e_1^s+f_0(e_2,e_3,e_1)e_2^s+f_0(e_3,e_1,e_2)e_3^s}
{4+f_2(e_1,e_2,e_3)e_1^s+f_2(e_2,e_3,e_1)e_2^s+f_2(e_3,e_1,e_2)e_3^s},\\
\hat{\pi}_2&=&\hat{\pi}_3=\frac{1+f_1(e_1,e_2,e_3)e_1^s+f_1(e_2,e_3,e_1)e_2^s+f_1(e_3,e_1,e_2)e_3^s}
{4+f_2(e_1,e_2,e_3)e_1^s+f_2(e_2,e_3,e_1)e_2^s+f_2(e_3,e_1,e_2)e_3^s},
\end{eqnarray*}
where
\begin{eqnarray*}
f_0(x,y,z)&:=&[1+\kappa-2\lambda-(1+\kappa)x]g(y,z)/[2(1+\kappa)^2\lambda(1+\lambda)(x-y)(x-z)], \\
f_1(x,y,z)&:=&[(1-\lambda)(1+\kappa-\lambda-\lambda^2)-(1+\lambda)(1-\kappa^2+\kappa\lambda-\lambda^2)x \\
&&\qquad{}+(1+\kappa)(1+\lambda)(\lambda-\kappa)x^2]g(y,z) \\
&&\quad{}/[2(1+\kappa)^2\lambda(1+\lambda)(1-\kappa\lambda)(x-y)(x-z)], \\
f_2(x,y,z)&:=&[2+2\kappa-4\lambda-2\kappa\lambda-\kappa^2\lambda+2\kappa\lambda^2+\lambda^3-(2+\kappa-\kappa^2+\lambda  \\
&&\qquad{}-2\kappa^2\lambda-\lambda^2+\kappa\lambda^2-\lambda^3)x+(1+\kappa)(1+\lambda)(\lambda-\kappa)x^2]g(y,z) \\
&&\quad{}/[(1+\kappa)^2\lambda(1+\lambda)(1-\kappa\lambda)(x-y)(x-z)],
\end{eqnarray*}
with
\begin{eqnarray*}
g(y,z)&:=&1+2\kappa-3\lambda+\kappa^2-2\kappa\lambda+\kappa^2\lambda-2\kappa\lambda^2+2\lambda^3 \\
& &\quad{}-(1+\kappa)(1+\lambda)(1+\kappa-2\lambda)(y+z)+(1+\kappa)^2(1+ \lambda)yz.
\end{eqnarray*}
Letting $\bm v:=(1/4)(1,\,1,\,-1,\,-1)$, we can write
$$ \hat{\bm \pi} = \bm u + (\hat{\bm \pi} - \bm u) = \bm u + ( 4 \hat{\pi}_0 -1) \bm v,$$
so $F_s = E_s + G_s H_s$, where
\begin{eqnarray*}
G_s&:=&4\hat{\pi}_0-1=2(\hat\pi_0-\hat\pi_2)\\
&=&2\frac{(f_0-f_1)(e_1,e_2,e_3)e_1^s+(f_0-f_1)(e_2,e_3,e_1)e_2^s+(f_0-f_1)(e_3,e_1,e_2)e_3^s}
{4+f_2(e_1, e_2, e_3)e_1^s+f_2(e_2,e_3,e_1)e_2^s+f_2(e_3, e_1, e_2) e_3^s}
\end{eqnarray*}
and
$$ H_s := {\bm v}{\bm R} \bm D_s{\bm L}{\bm\zeta}.$$
The argument that gave us (\ref{Es}) yields
\begin{eqnarray*}
H_{s} &=& b_0 - b_1 e_1^s - b_2 e_2^s - b_3 e_3^s,
\end{eqnarray*}
where
$$
b_0:=-\frac{1+\kappa-2\lambda-\lambda^2+\kappa\lambda^2}{4\lambda(1+\lambda)},
$$
and
$$
b_1:=h(e_1,e_2,e_3),\quad b_2:= h(e_2,e_3,e_1),\quad b_3:= h(e_3,e_1,e_2),
$$
with
\begin{eqnarray*}
h(x,y,z) &:=&[2+2\kappa-4\lambda-2\kappa\lambda-\kappa^2\lambda+2\kappa\lambda^2 +\lambda^3-(2+\kappa+\lambda-\kappa^2-\lambda^2 \\
&&\qquad{}-2\kappa^2\lambda+\kappa\lambda^2-\lambda^3)x+(1+\kappa)(1+\lambda)(\lambda-\kappa)x^2][1+3\kappa-2\lambda\\
&&\qquad{}+3\kappa^2-4\kappa\lambda-\lambda^2+\kappa^3-9\kappa\lambda^2+6\lambda^3+2\kappa^3\lambda-7\kappa^2\lambda^2+6\kappa\lambda^3\\
&&\qquad{}+\kappa^3\lambda^2-2\kappa^2\lambda^3+4\kappa\lambda^4-2\lambda^5-(1+\kappa)(1+\lambda)(1+2\kappa-3\lambda\\
&&\qquad{}+\kappa^2-2\kappa\lambda+\kappa^2\lambda-2\kappa\lambda^2+2\lambda^3)(y+z)\\ 
&&\qquad{}+(1+\kappa)^2(1+\lambda)^2(1+\kappa-2\lambda)yz]\\
&&\quad{}/[4(1+\kappa)^3\lambda(1+\lambda)^2(1-\kappa\lambda)(1-x)(x-y)(x-z)].
\end{eqnarray*}

Thus, (\ref{mu[r,s]}) and (\ref{mu[1,s]}) provide explicit, albeit complicated, formulas for $\mu_{[r,s]}(0)$ for all $r,s\ge1$.  They are valid provided only that $\beta^2+4\alpha^3\ne0$ (ensuring that $\bm P_B$ has distinct eigenvalues) and $\kappa\lambda\ne1$ (ensuring that the denominators of the formulas are nonzero).

We can extend the formulas to $\kappa\lambda=1$ by noting that, in this case, 0 is an eigenvalue of $\bm P_B$ and the two remaining nonunit eigenvalues, which can be obtained from the quadratic formula, are distinct from 0 and 1 unless $\kappa=\lambda=1$.  Here we are using the assumption that $\lambda<1+\kappa$, hence $\kappa>(-1+\sqrt{5})/2$.  This allows a spectral representation in such cases and again leads to formulas for $\mu_{[r,s]}(0)$ for all $r,s\ge1$, which we do not include here.  Writing the numerator of $f(x,y,z)$ temporarily as $(\lambda-\kappa)p(x)q(y,z)$, notice that the two nonunit, nonzero eigenvalues coincide, when $\kappa\lambda=1$, with the roots of $p(x)$, and this ordered pair of eigenvalues is also a zero of $q$.  This explains why the singularity on the curve $\kappa\lambda=1$ is removable. 

We cannot prove the analogue of Theorem 7 in this setting, so we state it as a conjecture. 

\begin{conjecture}
Let games $A$ and $B$ be as in Theorem 4 (with the bias parameter absent).  For every pair of positive integers $r$ and $s$, $\mu_{[r,s]}(0)>0$ if $\kappa<\lambda<1$ or $\kappa>\lambda>1$, $\mu_{[r,s]}(0)=0$ if $\kappa=\lambda$ or $\lambda=1$, and $\mu_{[r,s]}(0)<0$ if $\lambda<\min(\kappa,1)$ or $\lambda>\max(\kappa,1)$.  
\end{conjecture}

\begin{corollary}
Assume that the conjecture is true for some pair $(\kappa,\lambda)$ satisfying $0<\kappa<\lambda<1$ or $\kappa>\lambda>1$.
Let games $A$ and $B$ be as in Corollary 5 (with the bias parameter present).
For every pair of positive integers $r$ and $s$, there exists $\eps_0>0$, depending on $\kappa$, $\lambda$, $r$, and $s$, such that Parrondo's paradox holds for games $A$, $B$, and pattern $[r,s]$, that is, $\mu_A(\eps)<0$, $\mu_B(\eps)<0$, and $\mu_{[r,s]}(\eps)>0$, whenever $0<\eps<\eps_0$. 
\end{corollary}

We can prove a very small part of the conjecture.  Let $K$ be the set of positive fractions with one-digit numerators and denominators, that is, $K:=\{k/l: k,l=1,2,\ldots,9\}$, and note that $K$ has 55 distinct elements.

\begin{theorem}
The conjecture is true if \emph{(a)} $\kappa=\lambda>0$, \emph{(b)} $\kappa>0$ and $\lambda=1$, or \emph{(c)} $\kappa,\lambda\in K$, $\lambda<1+\kappa$, $\kappa\ne\lambda$, $\lambda\ne1$, $\beta^2+4\alpha^3\ne0$, and $\kappa\lambda\ne1$.
\end{theorem}

\begin{remark}
Parts (a) and (b) are the fair cases. Part (c) treats 2123 unfair cases (702 winning, 1421 losing), including the case $\kappa=1/9$ and $\lambda=1/3$ studied by Parrondo, Harmer, and Abbott (2000).  (The assumption that $\beta^2+4\alpha^3\ne0$ is redundant.)  The proof of part (c) is computer assisted.
\end{remark}

\begin{proof}
We begin with case (a), $\kappa=\lambda>0$.  In this case $c_0=c_1=c_2=c_3=0$, hence $E_s=0$, while
\begin{eqnarray*}
&&(f_0-f_1)(x,y,z)\\
&&\;\;{}=-{(\lambda-\kappa)[\lambda(\lambda-\kappa)-(1+\kappa-\lambda-\lambda^2)x+(1+\kappa)(1+\lambda)x^2]g(y,z)\over2(1+\kappa)^2\lambda(1+\lambda)(1-\kappa\lambda)(x-y)(x-z)},
\end{eqnarray*}
hence $G_s=0$.  We conclude that $\mu_{[r,s]}(0)=0$ for all $r,s\ge1$ in this case.  The argument fails only when $\kappa=\lambda=1$ because, assuming that $\kappa=\lambda$, only then is $\beta^2+4\alpha^3=0$ or is $\kappa\lambda=1$.  But in that case game $B$ is indistinguishable from game $A$, so again $\mu_{[r,s]}(0)=\mu_A(0)=0$ for all $r,s\ge1$.

Next we consider case (b), $\lambda=1$ and either $0<\kappa<1$ or $\kappa>1$.  We write the numerator of $f(x,y,z)$ temporarily as $(1-\kappa)p(x)q(y,z)$.
Each nonunit eigenvalue of $\bm P_B$ is a root of the cubic equation 
$x^3+a_2x^2+a_1x+a_0=0$ with coefficients (\ref{cubic-coeffs}).  With $\lambda=1$, this equation, multiplied by $4(1+\kappa)^2$, becomes
\begin{eqnarray*}
0&=&4(1+\kappa)^2x^3+4(1+\kappa)(1-\kappa)x^2+(3+\kappa)(1-\kappa)x+(1-\kappa)^2\\
&=&[2(1+\kappa)x^2+(1-\kappa)x+1-\kappa][2(1+\kappa)x+1-\kappa]\\
&=&p(x)[2(1+\kappa)x+1-\kappa].
\end{eqnarray*}
Moreover, we have $\beta^2+4\alpha^3=108(1+\kappa)^2(1-\kappa)^3(7+9\kappa)$.
If $\kappa<1$, then $\beta^2+4\alpha^3>0$ (the eigenvalues are distinct with $e_2$ and $e_3$ complex conjugates).  In this case, $p$ has complex roots, so $p(e_2)=p(e_3)=0$ and $e_1=-(1-\kappa)/[2(1+\kappa)]$; in addition, $e_2$ and $e_3$, being roots of $p$, have sum $-(1-\kappa)/[2(1+\kappa)]$ and product $(1-\kappa)/[2(1+\kappa)]$, hence $q(e_2,e_3)=0$.  Therefore, $c_0=c_1=c_2=c_3=0$ and $E_s=0$.
Finally, $(f_0-f_1)(x,y,z)$ has $p(x)$ as a factor and $g(e_2,e_3)=0$, hence $G_s=0$.
On the other hand, if $\kappa>1$, then $\beta^2+4\alpha^3<0$ (the eigenvalues are real and distinct).  In this case, two of the three eigenvalues $e_1$, $e_2$, and $e_3$ are roots of $p$, so the argument just given in the case $\kappa<1$ applies in this case as well with the two mentioned eigenvalues replacing $e_2$ and $e_3$.

We turn to case (c), but let us begin by treating the general case.  Given $\kappa>0$ and $\lambda>0$ with $\lambda<1+\kappa$, $\kappa\ne\lambda$, and $\lambda\ne1$, we clearly have $c_0\ne0$.  The conjecture says that $E_s$ and $F_s:=E_s+G_sH_s$ have the same sign as $c_0$ for all $s\ge1$.

We assume in addition that $\beta^2+4\alpha^3\ne0$ and $\kappa\lambda\ne1$.
If $\beta^2+4\alpha^3>0$, then we can take $e_1$ real and $e_2$ and $e_3$ complex conjugates.  It follows that $c_1$ is real and $c_2$ and $c_3$ are complex conjugates.  Similarly, $b_1$ is real and $b_2$ and $b_3$ are complex conjugates.  Therefore,
\begin{eqnarray*}
E_{s}&=&c_0-c_1 e_1^s-c_2 e_2^s-c_3 e_3^s=c_0-c_1 e_1^s-2\,{\rm Re}(c_2 e_2^s),\\
H_{s}&=&b_0-b_1 e_1^s-b_2 e_2^s-b_3 e_3^s=b_0-b_1 e_1^s-2\,{\rm Re}(b_2 e_2^s).
\end{eqnarray*}
If $\beta^2+4\alpha^3<0$, then $e_1$, $e_2$, and $e_3$ are real and distinct, as we have seen, and therefore $c_1,c_2,c_3,b_1,b_2,b_3$ are real.  In either case, $E_s$ has the same sign as $c_0$ if
\begin{equation}\label{EsBound}
|c_0|-\{|c_1|\,|e_1|^s+|c_2|\,|e_2|^s+|c_3|\,|e_3|^s\}
\end{equation}
is positive.
Notice that (\ref{EsBound}) increases in $s$ and is positive for $s$ large enough.  Given $(\kappa,\lambda)$ with $c_0\ne0$, $\beta^2+4\alpha^3\ne0$, and $\kappa\lambda\ne1$, there exists $s_0\ge1$, depending on $(\kappa,\lambda)$, such that (\ref{EsBound}) is positive for all $s\ge s_0$.  It is then enough to verify by direct computation that $E_s$ has the same sign as $c_0$ for $s=1,2,\ldots,s_0-1$. In fact, we believe, but cannot prove, that $s_0$ is never larger than 3.  If true, this would imply that the conjecture holds for $r\ge2$ and $s\ge1$.  (We are using the observations that $\mu_{[r,s]}(0)$ is continuous in $(\kappa,\lambda)$ and that the set of $(\kappa,\lambda)$ such that $\beta^2+4\alpha^3\ne0$ and $\kappa\lambda\ne1$ is dense.  We are also implicitly using the fact, not yet proved, that $\mu_{[r,s]}(0)\ne0$ on the curves $\beta^2+4\alpha^3=0$ and $\kappa\lambda=1$, except at $\kappa=\lambda=1$.)

The case $r=1$ is more complicated. Observe that $F_s$ has the same sign as $c_0$ if
\begin{eqnarray}\label{FsBound}
&&\!\!\!\!\!|c_0|-\bigg\{|c_1|\,|e_1|^s+|c_2|\,|e_2|^s+|c_3|\,|e_3|^s \nonumber\\
&&\!\!\!\!\!{}+2\frac{|f_0-f_1|(e_1,e_2,e_3)|e_1|^s+|f_0-f_1|(e_2,e_3,e_1)|e_2|^s+|f_0-f_1|(e_3,e_1,e_2)|e_3|^s}
{4-|f_2|(e_1, e_2, e_3)|e_1|^s-|f_2|(e_2,e_3,e_1)|e_2|^s-|f_2|(e_3, e_1, e_2)|e_3|^s} \nonumber\\
&&\qquad\qquad\qquad{}\cdot(|b_0|+|b_1|\,|e_1|^s+|b_2|\,|e_2|^s+|b_3|\,|e_3|^s)\bigg\} 
\end{eqnarray}
is positive, a stronger condition than requiring that (\ref{EsBound}) be positive.  Notice that (\ref{FsBound}) increases in $s$ and is positive for $s$ large enough.  Given $(\kappa,\lambda)$, again with $c_0\ne0$, $\beta^2+4\alpha^3\ne0$, and $\kappa\lambda\ne1$, there exists $s_1\ge1$, depending on $(\kappa,\lambda)$, such that (\ref{FsBound}) holds for all $s\ge s_1$.  It is then enough to verify by direct computation that $F_s$ has the same sign as $c_0$ for $s=1,2,\ldots,s_1-1$.  It appears that $s_1$ can be quite large.

We have carried out the required estimates in \textit{Mathematica} for the 2123 cases of part (c).  There were no exceptions to the conjecture.  Table 1 lists a few of these cases.
\end{proof}

\begin{table}
\caption{A few of the 2123 cases treated by Theorem 10(c).
$s_0$ is the smallest positive integer $s$ that makes (\ref{EsBound}) positive.
$s_1$ is the smallest positive integer $s$ that makes (\ref{FsBound}) positive.
\medskip}
\catcode`@=\active \def@{\hphantom{0}}
\renewcommand{\arraystretch}{1.}
\begin{center}
\begin{tabular}{ccccc}\hline
\noalign{\smallskip}
$(\kappa,\lambda)$ & region & $(p_0,p_1,p_3)$ & $s_0$ & $s_1$\\
& of Fig.\ 1 &&& \\
\noalign{\smallskip}\hline
\noalign{\smallskip}
$(1/9,1/3)$ & 1 & $(9/10,1/4,7/10)$ & 1 & 2 \\
$(1/3,1/9)$ & 6 & $(3/4,1/10,11/12)$ & 1 & 6 \\
$(9,3)$     & 3 & $(1/10,3/4,7/10)$ & 1 & 3 \\
$(1/9,1/8)$ & 1 & $(9/10,1/9,71/80)$ & 1 & 6 \\
$(1/9,8/9)$ & 1 & $(9/10,8/17,1/5)$ & 2 & 3 \\
$(8,1/9)$   & 6 & $(1/9,1/10,80/81)$ & 1 & 27 \\
$(4,9/2)$   & 2 & $(1/5,9/11,1/10)$ & 1 & 3 \\
$(3,3/2)$   & 4 & $(1/4,3/5,5/8)$ & 1 & 1 \\
$(3,2/3)$   & 5 & $(1/4,2/5,5/6)$ & 1 & 2 \\
\noalign{\smallskip}
\hline
\end{tabular}
\end{center}
\end{table}

Here are four of the simplest cases:
\begin{eqnarray*}
\mu_{[1,1]}(\eps)&=&{(\lambda-\kappa)(1-\lambda)\over2(2+\kappa+\lambda)(1+\lambda)}+O(\eps),\\
\mu_{[1,2]}(\eps)&=&{(\lambda-\kappa)(1-\lambda)(\kappa+\lambda+2\kappa\lambda)(2+2\kappa+\kappa\lambda-\lambda^2)\over3(1+\kappa)(1+\lambda)(2+\kappa+\lambda)(\kappa+3\lambda+4\kappa\lambda+\kappa\lambda^2-\lambda^3)}+O(\eps),\\
\mu_{[r,1]}(\eps)&=&{(\lambda-\kappa)(1-\lambda)\over2(r+1)(1+\kappa)(1+\lambda)}+O(\eps),\quad r\ge2,\\
\mu_{[r,2]}(\eps)&=&{(\lambda-\kappa)(1-\lambda)(1+2\kappa-\lambda)\over2(r+2)(1+\kappa)^2(1+\lambda)}+O(\eps),\quad r\ge2.
\end{eqnarray*}
Notice that the final factor in the numerator and the final factor in the denominator of $\mu_{[1,2]}(0)$ are both positive.  (Consider two cases, $\lambda\le\kappa$ and $\kappa<\lambda<1+\kappa$.)  The same is true of $\mu_{[r,2]}(0)$, $r\ge2$.

We turn finally to the evaluation of the asymptotic variance per game played of the player's cumulative profit.  We denote this variance by $\sigma^2_{[r,s]}(\eps)$, and we note that it suffices to assume that $\eps=0$ in the calculation to obtain the formula up to $O(\eps)$.

A formula for $\sigma_{[r,s]}^2(0)$ analogous to (\ref{mu[r,s]}) and (\ref{mu[1,s]}) would be extremely complicated.  We therefore consider the matrix formulas of Section 5 to be in final form.

For example, we have
\begin{eqnarray*}
\sigma^2_{[1,1]}(\eps)&=&(16+32\kappa+48\lambda+18\kappa^2+124\kappa\lambda+18\lambda^2+\kappa^3+91\kappa^2\lambda+71\kappa\lambda^2\\&&\quad{}-3\lambda^3+22\kappa^3\lambda+28\kappa^2\lambda^2+38\kappa\lambda^3-8\lambda^4+\kappa^3\lambda^2+15\kappa^2\lambda^3\\
&&\quad{}-\kappa\lambda^4+\lambda^5)/[2(1+\lambda)^2(2+\kappa+\lambda)^3]+O(\eps),\\
\sigma^2_{[r,1]}(\eps)&=&1-\frac{(\lambda-\kappa)(8+7\kappa+17\lambda+18\kappa\lambda+6\lambda^2+7\kappa\lambda^2+\lambda^3)}
{4(r+1)(1+ \kappa)^2(1+\lambda)^2}+O(\eps),\\
\sigma^2_{[r,2]}(\eps)&=&1+(8+36\kappa-20\lambda+71\kappa^2-42\kappa\lambda-37\lambda^2+64\kappa^3+24\kappa^2\lambda\\
&&\qquad{}-120\kappa\lambda^2+20\kappa^4+96\kappa^3\lambda-118\kappa^2\lambda^2-12\kappa\lambda^3+6\lambda^4\\
&&\qquad{}+48\kappa^4\lambda-16\kappa^3\lambda^2-24\kappa^2\lambda^3+4\kappa\lambda^4+4\lambda^5+20\kappa^4\lambda^2\\
&&\qquad{}-16\kappa^3\lambda^3-k^2\lambda^4+6\kappa\lambda^5-\lambda^6)/[4(r+2)(1+\kappa)^4(1+\lambda)^2]\\
&&\;\;{}+O(\eps),
\end{eqnarray*}
for $r\ge2$.  The formula for $\sigma^2_{[1,2]}(\eps)$ is omitted above; $\sigma^2_{[1,2]}(0)$ is the ratio of two polynomials in $\kappa$ and $\lambda$ of degree 16, the numerator having 110 terms.

When $\kappa = 1/9$ and $\lambda = 1/3$, the mean and variance formulas simplify to
\begin{align*}
\mu_{[1,1]}(\eps)&={1\over44}+O(\eps),&\sigma^2_{[1,1]}(\eps)&={8945\over10648}+O(\eps),\\
\mu_{[1,2]}(\eps)&={203\over16500}+O(\eps),&\sigma^2_{[1,2]}(\eps)&={1003207373\over998250000}+O(\eps),\\
\mu_{[2,1]}(\eps)&={1\over60}+O(\eps),&\sigma^2_{[2,1]}(\eps)&={1039\over1200}+O(\eps),\\
\mu_{[2,2]}(\eps)&={1\over100}+O(\eps),&\sigma^2_{[2,2]}(\eps)&={19617\over20000}+O(\eps).
\end{align*}

\section{Why does Parrondo's paradox hold?}

Several authors (e.g., Ekhad and Zeilberger 2000, Rahmann 2002, Philips and Feldman 2004) have questioned whether Parrondo's paradox should be called a paradox.  Ultimately, this depends on one's definition of the word ``paradox,'' but conventional usage (e.g., the St.\ Petersburg paradox, Bertrand's paradox, Simpson's paradox, and the waiting-time paradox) requires not that it be unexplainable, only that it be surprising or counterintuitive. Parrondo's paradox would seem to meet this criterion.

A more important issue, addressed by various authors, is, Why does Parrondo's paradox hold?  Here we should distinguish between the random mixture version and the nonrandom pattern version of the paradox.  The former is easy to understand, based on an observation of Moraal (2000), subsequently elaborated by Costa, Fackrell, and Taylor (2005).  See also Behrends (2002).  Consider the capital-dependent games first.  The mapping $f(\rho):=(\rho^2/(1+\rho^2),1/(1+\rho))=(p_0,p_1)$ from $(0,\infty)$ into the unit square $(0,1)\times(0,1)$ defines a curve representing the fair games, with the losing games below and the winning games above.  The interior of the line segment from $f(1)=({1\over2},{1\over2})$ to $f(\rho)$ lies in the winning region if $\rho<1$ (i.e., $p_0<{1\over2}$) and in the losing region if $\rho>1$ (i.e., $p_0>{1\over2}$), as can be seen by plotting the curve (as in, e.g., Harmer and Abbott 2002).  These line segments correspond to the random mixtures of game $A$ and game $B$.  Actually, it is not necessary to plot the curve.  Using (\ref{profit-param}), we note that the curve defined by $f$ is simply the graph of $g(p_0):=1/(1+\sqrt{p_0/(1-p_0)}\,)$ ($0<p_0<1$), and by calculating $g''$ we can check that $g$ is strictly convex on $(0,{1\over2}]$ and strictly concave on $[{1\over2},1)$.

The same kind of reasoning applies to the history-dependent games, except that now the mapping $f(\kappa,\lambda):=(1/(1+\kappa),\lambda/(1+\lambda),1-\lambda/(1+\kappa))=(p_0,p_1,p_3)$ from $\{(\kappa,\lambda):\kappa>0,\lambda>0,\lambda<1+\kappa\}$ into the unit cube $(0,1)\times(0,1)\times(0,1)$ defines a surface representing the fair games, with the losing games below and the winning games above.  The interior of the line segment from $f(1,1)=({1\over2},{1\over2},{1\over2})$ to $f(\kappa,\lambda)$ lies in the winning region if $\kappa<\lambda<1$ (i.e., $p_0+p_1>1$ and $p_1<{1\over2}$) or $\kappa>\lambda>1$ (i.e., $p_0+p_1<1$ and $p_1>{1\over2}$).  It is possible but difficult to see this by plotting the surface on a computer screen using 3D graphics.  Instead, we note using (\ref{history-param}) that the surface defined by $f$ is simply the graph of $g(p_0,p_1):=1-p_0p_1/(1-p_1)$ ($0<p_0<1$, $0<p_1<1$, $p_0p_1<1-p_1$).  
So with $h(t):=g((1-t)/2+tp_0,(1-t)/2+tp_1)$ for $0\le t\le 1$, we compute 
$$
h''(t)=-{4(p_0+p_1-1)(2p_1-1)\over[1-(2p_1-1)t]^3}
$$
and observe that $h''(t)>0$ for $0\le t\le 1$ if and only if both $p_0+p_1>1$ and $p_1<{1\over2}$, or both $p_0+p_1<1$ and $p_1>{1\over2}$.  In other words, $g$ restricted to the line segment from $({1\over2},{1\over2})$ to $(p_0,p_1)$ is strictly convex under these conditions.  This confirms the stated assertion.

The explanations for the nonrandom pattern version of the paradox are less satisfactory.  Ekhad and Zeilberger (2000) argued that it is because ``matrices (usually) do not commute.''  The ``Boston interpretation'' of H. Eugene Stanley's group at Boston University claims that it is due to noise.  As Kay and Johnson (2003) put it, ``losing cycles in game B are effectively broken up by the memoryless behavior, or `noise', of game A.''  Let us look at this in more detail.  

We borrow some assumptions and notation from the end of Section 5.  Assume that $\eps=0$.  Denote by $\bm P_A'$ the matrix with $(i,j)$th entry $(\bm P_A)_{ij}w(i,j)$, and assume that $\bm P_A'$ has row sums equal to 0.  Denote by $\bm P_B'$ the matrix with $(i,j)$th entry $(\bm P_B)_{ij}w(i,j)$, and define $\bm \zeta:=\bm P_B'(1,1,\ldots,1)^\T$ to be the vector of row sums of $\bm P_B'$.  Let $\bm\pi_B$ be the unique stationary distribution of $\bm P_B$, and let $\bm \pi_{s,r}$ be the unique stationary distribution of $\bm P_B^s\bm P_A^r$.  Then the asymptotic mean per game played of the player's cumulative profit for the pattern $[r,s]$ when $\eps=0$ is
$$
\mu_{[r,s]}(0)=(r+s)^{-1}\bm \pi_{s,r}(\bm I+\bm P_B+\cdots+\bm P_B^{s-1})\bm\zeta.
$$
If it were the case that $\bm\pi_{s,r}=\bm\pi_B$, then we would have 
$$
\mu_{[r,s]}(0)=(r+s)^{-1}s\bm\pi_B\bm\zeta=0,
$$
and the Parrondo effect would not appear.  If we attribute the fact that
$\bm\pi_{s,r}$ is not equal to $\bm\pi_B$ to the ``noise'' caused by game $A$, then we can justify the Boston interpretation.

The reason this explanation is less satisfactory is that it gives no clue as to the sign of $\mu_{[r,s]}(0)$, which indicates whether the pattern $[r,s]$ is winning or losing.  We propose an alternative explanation (the Utah interpretation?) that tends to support the Boston interpretation.  To motivate it, we observe that $(r+s)\mu_{[r,s]}(0)$ can be interpreted as the asymptotic mean \textit{per cycle of $r$ plays of game $A$ and $s$ plays of game $B$} of the player's cumulative profit when $\eps=0$.  If $s$ is large relative to $r$, then the $r$ plays of game $A$ might reasonably be interpreted as periodic noise in an otherwise uninterrupted sequence of plays of game $B$.  We will show that
$$
\lim_{s\to\infty}(r+s)\mu_{[r,s]}(0)
$$
exists and is finite for every $r\ge1$.  If the limit is positive for some $r\ge1$, then $\mu_{[r,s]}(0)>0$ for that $r$ and all $s$ sufficiently large.  If the limit is negative for some $r\ge1$, then $\mu_{[r,s]}(0)<0$ for that $r$ and all $s$ sufficiently large.  These conclusions are weaker than those of Theorem 7 and the conjecture but the derivation is much simpler, depending only on the fundamental matrix of $\bm P_B$, not on its spectral representation.

Here is the derivation.  First, $\bm\pi_{s,r}=\bm\pi\bm P_A^r$, where $\bm\pi$ is the unique stationary distribution of $\bm P_A^r\bm P_B^s$.  Now
$$
\lim_{s\to\infty}\bm P_A^r\bm P_B^s=\bm P_A^r\bm\Pi_B=\bm\Pi_B,
$$
where $\bm\Pi_B$ is the square matrix each of whose rows is $\bm\pi_B$.  We conclude that $\bm \pi_{s,r}\to \bm\pi_B\bm P_A^r$ as $s\to\infty$ and therefore that
\begin{eqnarray*}
\lim_{s\to\infty}(r+s)\mu_{[r,s]}(0)&=&\lim_{s\to\infty}\bm \pi_{s,r}(\bm I+\bm P_B+\cdots+\bm P_B^{s-1})\bm\zeta\\
&=&\bm\pi_B\bm P_A^r\lim_{s\to\infty}\bigg[\sum_{n=0}^{s-1}(\bm P_B^n-\bm\Pi_B)+s\bm\Pi_B\bigg]\bm\zeta\\
&=&\bm\pi_B\bm P_A^r\bigg[\sum_{n=0}^\infty(\bm P_B^n-\bm\Pi_B)\bigg]\bm \zeta\\
&=&\bm\pi_B\bm P_A^r(\bm Z_B-\bm \Pi_B)\bm \zeta\\
&=&\bm\pi_B\bm P_A^r\bm Z_B\bm \zeta,\qquad r\ge1,
\end{eqnarray*}
where $\bm Z_B$ denotes the fundamental matrix of $\bm P_B$ (see (\ref{Z})).  Since $\bm\pi_B\bm Z_B=\bm\pi_B$ and $\bm\pi_B\bm\zeta=0$, it is the $\bm P_A^r$ factor, or the ``noise'' caused by game $A$, that leads to a (typically) nonzero limit.

In the capital-dependent setting we can evaluate this limit as
$$
\lim_{s\to\infty}(r+s)\mu_{[r,s]}(0)={3a_r(1-\rho)^3(1+\rho)\over2(1+\rho+\rho^2)^2},\qquad r\ge1,
$$
where $a_r$ is as in (\ref{a_r}), while in the history-dependent setting it becomes
$$
\lim_{s\to\infty}(r+s)\mu_{[r,s]}(0)=\frac{(1+\kappa)(\lambda-\kappa)(1-\lambda)}{4\lambda(2+\kappa+\lambda)},\qquad r\ge1.
$$
Of course we could derive these limits from the formulas in Sections 6 and 7, but the point is that they do not require the spectral representation of $\bm P_B$---they are simpler than that.  They also explain why the conditions on $\rho$ are the same in Theorems 2 and 7,  and the conditions on $\kappa$ and $\lambda$ are the same in Theorem 4 and the conjecture.

\section*{Acknowledgments}

The research for this paper was carried out during J. Lee's visit to the Department of Mathematics at the University of Utah in 2008--2009.  The authors are grateful to the referees for valuable suggestions.

\end{document}